\documentclass[dvips]{article}
\usepackage{graphicx}

\font\BBb=msbm10 at 12pt

\newcommand{\square}{\hbox{${\vcenter{\hrule height.4pt
  \hbox{\vrule width.4pt height6pt \kern6pt
     \vrule width.4pt}
  \hrule height.4pt}}$}}

\newcommand {\qed} {\hfill \nobreak \square \medbreak}

\newcommand{\be}{\begin{equation}}
      \newcommand{\ee}{\end{equation}}
      \newcommand{\ba}{\begin{eqnarray}}
       \newcommand{\ea}{\end{eqnarray}}
\newcommand{\ban}{\begin{eqnarray*}}
\newcommand{\ean}{\end{eqnarray*}}
\newcommand{\Bbb}[1]{\mbox{\BBb #1}}

\newcommand{\ul}{\underline}

\newcommand{\NN}{\Bbb{N}}
\newcommand{\RR} {\Bbb{R}}
\newcommand{\TT} {\Bbb{T}}
\newcommand{\SSS} {\Bbb{S}}
\newcommand{\HH} {\Bbb{H}}
\newcommand{\EE} {\Bbb{E}}

 \newcommand{\Pf}{\noindent {\bf Proof:} }

\newcommand{\sect}[1]{\section{#1} \setcounter{equation}{0}}
\newtheorem{theorem}{Theorem}[section]
\newtheorem{theo}{Theorem}[section]
\newtheorem{remark}[theo]{Remark}
\newtheorem{defn}{Definition}[section]
\newtheorem{example}{Example}[section]

\newtheorem{lemma}[theo]{Lemma}

\newtheorem{coro}[theo]{Corollary}
\newtheorem{corollary}[theo]{Corollary}

\begin{document}
\title{Friedmann Cosmology and Almost Isotropy}

\author{Christina Sormani\thanks {Partially supported by NSF Grant \# DMS-0102279
and a PSC CUNY 34 Award}}

\date{}
\maketitle

\noindent {\bf Abstract:}
In the Friedmann Model of the universe, cosmologists assume that
spacelike slices of the universe are Riemannian manifolds of constant sectional curvature.
This assumption is justified via Schur's Theorem by stating that the spacelike
universe is locally isotropic. Here we define  a Riemannian
manifold as almost locally isotropic in a sense which allows both 
weak gravitational lensing in all directions and strong gravitational
lensing in localized angular regions at most points.  We then prove that such a manifold is Gromov-Hausdorff
close to a length space $Y$ which is a collection of space forms joined at discrete points.
Within the paper we define a concept we call an ``exponential length space''
and prove that if such a space is locally isotropic then it is a space form.

\sect{Introduction} \label{intro}

The Friedmann Model of the universe is a Lorentzian manifold satisfying
Einstein's equations which is assumed to be a warped product of a
a space form 
with the real line.
\cite{Fra} \cite{Peeb} \cite{CW}.
Recall that a space form is a complete Riemannian manifold with constant
sectional curvature.
This assumption is ``justified'' in the references above by stating that the spacelike
universe is locally isotropic:

\begin{defn} \label{iso} {\em
A Riemannian manifold $M$ is $R$ {\em locally isotropic} if for all $p\in M$
and  for every element $g\in SO(n,\RR)$
there is an isometry between balls, $f_g:B_p(R)\to B_p(R)$, such that
$f_g(p)=p$ and $df_g=g:TM_p\to TM_p$. }
\end{defn}

Clearly if $M$ is locally isotropic
then its sectional curvature $K_p(\sigma)$
depends only on $p$ not the $2$ plane $\sigma \subset TM_p$
and, by Schur's lemma, it has constant sectional curvature
(c.f. \cite[p.142] {Fra}).  

Now the universe is not exactly locally
isotropic and is only an approximately so.
To deal with this, cosmologists test perturbations
of the Friedmann model and look for measurable effects
on light rays.  The most popular perturbation is the
Swiss Cheese model in which holes are cut out of
the standard model and replaced with Schwarzschild solutions
\cite{Kan}\cite{DyRo}.  The effects of these clumps of mass
have been tested using random distribution
\cite{HoWa} and fractal distribution \cite{GabLab} of the massive regions.  However
all these studies of possible cosmologies are making the
assumption that the Friedmann model is stable in some sense.

It should be noted that  Schur's Lemma is not stable.  Noncompact examples
by Gribkov and compact examples by Currier 
show that Riemannian manifolds 
whose sectional
curvature satisfies
\be
|K_p(\sigma)-K_p|<\epsilon \qquad \forall \textrm{ 2 planes }\sigma \subset TM_p
\ee
can still have 
\be
\max_{p\in M}|K_p| - \min_{q\in M} |K_q| =1,
\ee
and thus do not have almost constant sectional curvature
\cite{Grib} \cite{Cur}.
The only stability theorem for Schur's Lemma has been
proven by Nikolaev, and it makes an integral approximation
on the pointwise sectional curvature variation \cite{Nik}.
Furthermore before one could even apply Nikolaev's 
Stability Theorem in our situation,
one would need to investigate whether 
a space which is almost isotropic in some sense has almost constant
sectional curvature at each point.  

In this paper we show that the implication that a locally isotropic
Riemannian manifold is a space form is stable with respect to the Gromov Hausdorff
topology on Riemannian manifolds.  This stability uses a definition of
almost isotropy which the author has constructed to allow
observed inhomogeneities in the universe including strong and weak
gravitational lensing as long as the weak lensing is very weak
and the strong lensing is localized in an angular sense.
[Theorem~\ref{MainThm} and Theorem~\ref{RicciThm}].
Furthermore the Swiss Cheese models of the universe
are almost isotropic in the sense that will be used
in this paper \cite{Kan}\cite{DyRo}.  

We will begin by providing a more angular rephrasing of the
definition of local isotropy
which is equivalent to Definition~\ref{iso}.

\begin{defn}\label{iso2} {\em
A Riemannian manifold $M$ is $R$ {\em locally isotropic} if for all $p\in M$
there is a radius $R_p>0$ less than the injectivity radius at $p$ and 
a function $F_p:[0,\pi]\times [0, R) \times [0,R)$ 
such that
for all unit vectors $v,w\in S^{n-1}\subset TM_p$ and for all $s,t \in (0,R_p)$ we have 
\be
d_M(exp_p(sv), exp_p(tw))=F_p(d_{S}(v, w), s, t).
\ee
where $d_S(v,w)$ is the angle between $v$ and $w$.
We will call $F_p$ the {\em isotropy function}
about $p$ and $R_p$ the {\em radius of isotropy at $p$}.
Furthermore $M$ is {\em uniformly isotropic}
on a region $U$, if $F_p$ is constant for $p\in U$.}
\end {defn}

Note that as described above, a locally isotropic manifold is a space form.
Thus $F_p$ must increase in its first variable and $F_p(\pi,t,t)=2t$
for $t$ sufficiently small.   While $F_p$ is not
assumed to be constant here, by Schur's Lemma, it must be constant
and in fact
\be
F_p(\theta, s,t)=F_K(\theta,s,t)
\ee
where $F_K(\theta,s,t)$ is the length of the third side 
of a triangle with angle $\theta$ between sides of lengths
$s$ and $t$ in $\HH^n$, $\SSS^n$ or $\EE^n$ of constant
sectional curvature $K$.  This is a well known function,
e.g.
\be
F_0(\theta,s,t)=s^2+t^2-2st \cos\theta.
\ee
 
We have now defined local isotropy as a property of geodesics
emanating from a point, a property that can be measured astronomically
if one assumes that light travels along spacial geodesics.  This is true
if we have given the spacelike slice of the universe the Fermat metric
(which incidentally is proportional to the restricted metric in the 
Friedmann model) c.f. pages 90-92 and 141-143 in \cite{Fra}.

In the following definition we approximate local stability in a sense
which will allow weak gravitational lensing of some geodesics and
strong gravitational lensing of those that enter a region $W \subset M$.  
Since geodesics entering $W$ behave unpredictably, we will restrict
our relatively good behavior to geodesics emanating from $p$ outside
a tubular neighborhood $T_\epsilon(W)$ of $W$.

\begin{figure}[htbp]
\includegraphics[height=2.5in ]{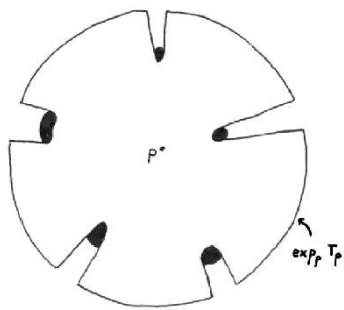}
\caption{ We've filled $W\subset M^2$  in black and outlined $exp_p(T_p)$ . } \label{fig1}
\end{figure}

\begin{defn} \label{almiso}  {\em
Given $\epsilon<1$, $R>1$, a Riemannian manifold $M^n$ and a subset $W\subset M^n$,
we say that $M^n$ is {\em locally $(\epsilon, R)$-almost isotropic off of}
$W$
if for all $p\in M\setminus T_\epsilon(W)$, we have a set of tangent vectors
\be \label{TpW}
T_p=T_{p,W}=\{ v\in B_0(R)\subset TM_p: exp_p([0,1]v)\cap W=\emptyset \}
\ee
and a function $F_p:[0,\pi]\times [0,R)\times [0,R) \to [0,2R)$
satisfying
\begin{eqnarray} 
(a) & F_p(\theta_1, a, b)<F_p(\theta_2, a,b) \quad \forall \theta_1<\theta_2 
\label{abc1a}\\
(b) & F_p(\pi, t,t) >t \label{abc1b}\\
(c) & F_p(0,0,R)=R \label{abc1c}
\end{eqnarray}
such that
\be \label{Fpdef}
|d_M(exp_p(v), exp_p(w))-F_p(d_S(v/|v|,w/|w|), |v|, |w|)| <\epsilon
\qquad \forall v,w\in T_p.
\ee

We will call $F_p$ the {\em almost isotropy function} about $p$ 
and $R$ the {\em isotropy radius}.  We will say that $M^n$ is {\em
uniformly almost isotropic} on a region $U$ if $F_p$ can be taken
to be constant for $p\in U$. 
See Figure~\ref{fig1}.}
\end{defn}

%\psdraft
%\psfull

% nocol is the only place we use $F_p(-1, t,t)>t $ really only need limit>0
%which may be a consequence of other parts of the defn if 
%we assume the limit isn't a single point 
% F(0,0,R)=R gives min geodesics for exp_y.

This definition captures the concept that the universe looks almost
the same in many directions as an angular view, $T_p$, but allows for
some directions to be poorly behaved after they pass through a
region, $W$, with strong gravitational lensing effects. Small gravitational lensing
is absorbed in the flexibility of (\ref{Fpdef}). 
Note that assumptions (a) (b) and (c) on $F$ all hold on isotropy functions.
Condition (b) guarantees that there is a geodesic
$exp_p(tv)$ whose end points are at least a distance $R$ apart.
This condition will replace the standard injectivity radius condition
often imposed on Riemannian manifolds when studying their limits.  By only requiring
that $d_M(exp_p(-Rv), exp_p(Rv))>R$ and not $=2R$, we are not demanding
that $p$ be a midpoint of any long minimizing geodesic as is the case
with an injectivity radius bound.  This is not a strong assumption to make
for certainly it would seem that there should be two opposing
directions in the sky that are far apart from each other.

We now wish to impose some restrictions on the size of the set $W$ where
$M$ fails to be almost isotropic using an angular
measurement.  The idea we are trying to capture is that very few directions
in the sky exhibit strong gravitational lensing.  

\begin{figure}[htbp]
\includegraphics[height=2.5in ]{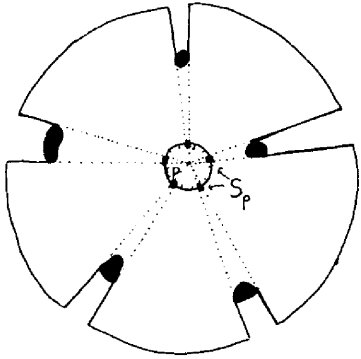}
\caption{We've drawn  $W\subset M^2$ in black, $S_p\subset S^1\subset TM_p$
marked as five intervals on a circle.} \label{fig2}
\end{figure}

\begin{defn} \label{unseen} 
{\em
A subset $W$ of $M$ is {\em $(\epsilon, R)$-almost unseen} if
for all $p\in M\setminus T_\epsilon(W)$, the set of directions
from $p$ passing through $W\cap B_p(R)$,
\be
S_p=\{v/|v|: v \in B_0(R)\setminus T_{p,W}\subset TM_p\},
\ee
where $T_{p,W}$ was defined in (\ref{TpW}),
is contained in a disjoint set of balls,
\be
S_p\subset \bigcup_{j=1}^{N} B_{w_j}(\epsilon_j) \subset S^{n-1}\subset TM_p
\ee
where
$
B_{w_j}(3\epsilon_j)
$
are disjoint and $\epsilon_j<\epsilon$. 
See Figure~\ref{fig2}.} %(figure 1 with $S_p$ emphasized)
\end{defn}

%\psdraft
%\psfull

In our theorems we will also assume that $W\subset\bigcup_k B_{q_k}(\epsilon)$
where $d(q_k, q_j)>2R$.
In some sense this means we are making the assumption that the ``black holes'' are
small (which they appear to be from a cosmic perspective) and are far between.
Ones which are closer together can be fit in a common $\epsilon$ ball.
Spaces with thin wormholes that are long do not satisfy this condition.
To allow for such spaces we can cut off the worm holes and smooth them out.
In which case we are really only concerned with the universe on ``our side''
of the wormholes, our connected region.  

We may now state the main theorem and then its cosmological implications.  Here 
$d_{GH}(X_1, X_2)$ denotes the
Gromov-Hausdorff distance between $X_1$ and $X_2$.  Section~\ref{sectlimits}
contains a description of this distance between spaces.
We will also use the notation $Ricci(M^n)\ge (n-1)H$ to denote that
the standard assumption that Ricci curvature is bounded below in the sense that
\be
Ricci_p(v,v) \ge (n-1)Hg(v,v) \qquad \forall p\in M^n, \,\,\forall v\in TM_p.
\ee 

\begin{theorem} \label{RicciThm}  
Given $H>0$, $n \in \NN$, $\bar{R},R>0$, $D>0$ and $\delta>0$ there exists 
\be
\epsilon=\epsilon(H,n, \bar{R},R, D, \delta)>0
\ee 
such that if $\bar{B}_p(D) \subset M^n$ is a closed ball in a complete
Riemannian manifold with the $Ricci(M^n) \ge (n-1)H$ such that
$M^n$ is locally $(\epsilon,R)$-almost isotropic off of $W$ where
$W$ is an $(\epsilon, R)$-almost unseen set
contained in uniformly disjoint balls,
\be \label{Wballs}
W\subset \bigcup_j B_{q_j}(\epsilon) \textrm{ where } B_{q_j}(\bar{R}) 
\textrm{ are disjoint,}
\ee 
then 
\be \label{ricciGH}
d_{GH}(B_p(D) \subset M^n, B_y(D)\subset Y)<\delta
\ee
where $Y$ is an $n$ dimensional Riemannian manifold with constant sectional
curvature $K\ge H$ and injectivity radius greater than $R$.  
Furthermore $M^n$ is uniformly $(\epsilon+\delta)$-almost 
isotropic on $B_p(D)$,
\be \label{ricciexp}
|F_q(\theta, s,t)-F_K(\theta,s,t)|<\delta \qquad \forall \, q\in B_p(D),   
\ee
where $F_K$ is the isotropy function of a simply connected
space form of sectional curvature $K$.
\end{theorem}

%\psdraft
\begin{figure}[htbp]
\includegraphics[height=2.2in ]{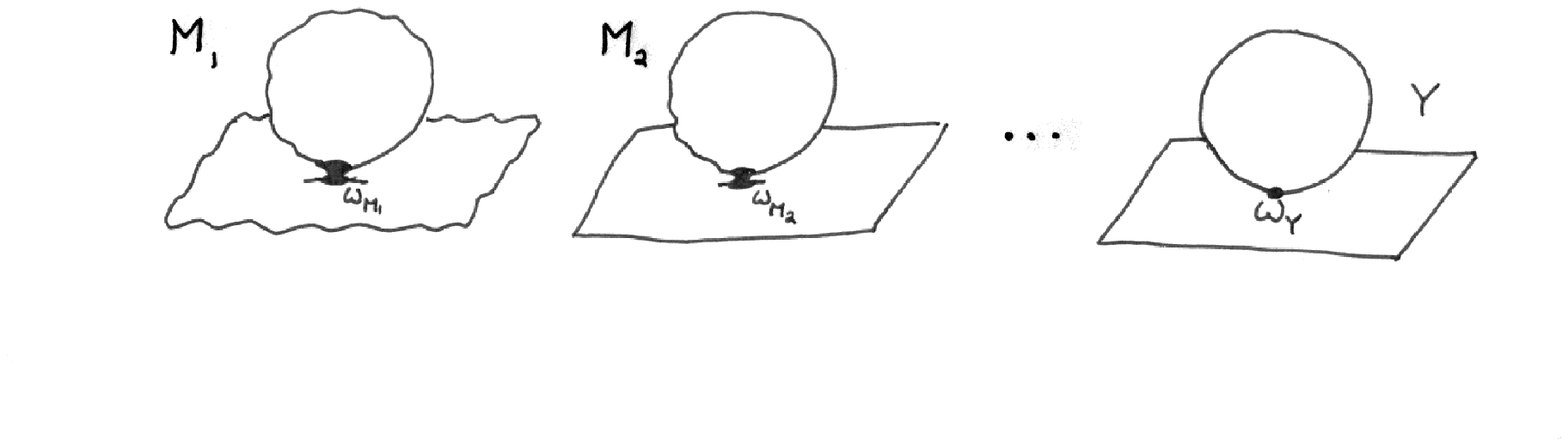}
\caption{Here we see Weak Gravitational Lensing all over the bumpy 
almost isotropic $M_i$.  There is Strong Gravitational 
Lensing at the necks of the $M_i$ which is restricted to the the subsets, $W_{M_i}$, marked in black.  
These $M_i$ converge in the Gromov Hausdorff sense to $Y$ with a single point in $W_Y$.  } 
\label{fig5}
\end{figure}
%\psfull

Note that if $M$ is compact then we can take $D$ to be the diameter
and we need not deal with the closed balls.  In this case $Y$ will be
a compact space form.  This does not require that $Y$ has positive curvature
since it may be a torus or a compact quotient of hyperbolic space.
If one further adds the condition that $M$ is a simply connected
compact manifold, then $Y$ must be a sphere (see Remark~\ref{soweirmk}).

If $M$ is noncompact,  $F_q$ may change slowly from point to 
point in $M$ such that the $K$ in (\ref{ricciexp}) depends on the ball.  See 
Example~\ref{ballchange}.    Nevertheless one can use (\ref{ricciexp})
to control the growth of the change in $K$ because the same $\delta$
holds for all balls.

Cosmologically, Theorem~\ref{RicciThm} says that if a space has sufficiently
small weak gravitational lensing and sufficiently localized
strong gravitational lensing
as viewed from most points in space and if the strong gravitational lensing
is caused by regions which are contained in sufficiently small balls
then in fact one can estimate the distances between stars 
whose light has not passed through regions of strong gravitational lensing
using the standard formulas involving only the angle between them as viewed 
from earth and the distance to the two stars (\ref{ricciexp}).  Of course
one must estimate $K$ as usual, but this can be done using astronomical
data measured from earth, and then the same 
$K$ can be used in all
directions and from any basepoint (not just earth).  

The first equation
(\ref{ricciGH}) is a bit more complicated to describe quickly
other than to say in some rough, not smooth sense the space is close
to a space of constant curvature.  A discussion of the Gromov-Hausdorff distance,
denoted $d_{GH}$ can be found in Section~\ref{sectlimits} and a good
reference is \cite{BBI}.

The following immediate corollary of Anderson's Smooth
Convergence Theorem clarifies this closeness if one adds
an additional upper Ricci curvature bound \cite{And}.

\begin{coro} \label{And}
Given $H>0$, $n, \in \NN$, $R>0$, $D>0$ $\delta>0$ there exists 
$\epsilon=\epsilon(H,n, R, D, \delta)>0$ such that if $B_p(D) \subset M^n$ is 
closed ball in a complete
Einstein Riemannian manifold with the $|Ricci| \le H$  $injrad\ge i_0$ $diam(M^n)\le D$
and
$M^n$ is locally $(\epsilon,R)$-almost isotropic,
then 
$B_p(D)$ is $C^\infty$ close to a ball in a compact space form.
\end{coro}

This $C^\infty$ closeness allows one to study the properties of the
universe using a smooth variation of the standard Friedmann model.
That is, spaces which are $C^\infty$ close to a space form can be studied
just by smoothly varying the metric on a space form and do not have possible
additional topology as might occur in the case when it is only Gromov-Hausdorff
close to a space form.

In some sense these two results are not as natural as one would hope because
one would like to permit a sequence of Schwarzschild universes that 
do not approach a Riemannian manifold but rather a pair of planes joined at a point.  Here we refer to a Riemannian Schwarzschild universe which is a
spacelike slice with a 0 second fundamental form in a spherically symmetric 
static Lorentzian Schwarzschild universe.  It is the universe which
achieves equality in the Penrose inequality \cite{SchYau}, and has one black 
hole with two asymptotically Euclidean ends (or two distinct `cosmos')
on either side.  If one
examines a sequence of such universes rescaled so that the black holes 
shrink to a point, the limit is a pair of planes joined at that point.
Such sequences of universes
do not have uniformly bounded Ricci curvature inside their 'black holes'.  

To apply the
above results to spaces which contain spherical black holes, like the 
Schwarzschild universes and the ones in 
the Swiss Cheese models,
we can take our manifold to be an edited
version of actual space in which all black holes have been cut out and replaced
with Euclidean balls that have been smoothly attached with bounded Ricci
curvature.  These black holes would then
be included in the region $W$ with strong
gravitational lensing so it does not matter that geodesics passing through
them no longer behave the way they did before the editing process.  The
new manifold would still be almost isotropic and then Theorem~\ref{RicciThm}
and Corollary~\ref{And} can be applied.

In the next theorem, we allow for an almost isotropic space which includes 
black
holes and any other sort of region that is badly behaved without having to cut
and paste the manifold.  This has the advantage that one need not make any assumption
on the local curvature of the space in these bad regions allowing for undiscovered
phenomenon like short wormholes or networks of wormholes or any other distortion
of space that is restricted to a collection of small balls.  It also matches
the conditions of the Swiss Cheese Models
of the universe studied in \cite{Kan}\cite{DyRo}.  

We replace the Ricci bound in Theorem~\ref{RicciThm} 
by a significantly more general ball packing assumption and obtain a 
slightly weaker result 
that neatly matches the idea of black holes joining pairs of
universes.

\begin{defn}
Given a Riemannian manifold and a map $f:(0,R)\times(0,\infty)\to \NN$
we say that $M$ has the $f$ {\em ball packing property} if
for any $s\in (0,R)$ and $t\in (0,\infty)$ the maximum number
of disjoint balls of radius $s$ contained in a ball of radius $t$
is bounded by $f(s,t)$.
\end{defn}

Note that Gromov's compactness theorem says that a sequence of
Riemannian manifolds has a converging subsequence iff there exists a function
$f$ such that all manifolds satisfy the same $f$ ball packing
property.  When a Riemannian manifold has $Ricci \ge -(n-1)H$ then by the Bishop-Gromov
Volume Comparison Theorem, 
$f$ is an explicit function of $H$ involving $sinh$ \cite{Bi}\cite{Gr}.
The limits of sequences of manifolds which satisfy the same $f$ ball packing
property are not necessarily manifolds but are complete length spaces (see
Definition~\ref{cls}).

\begin{theorem} \label{MainThm}  
Given  $n \in \NN$, $\bar{R}, R>0$, $D>0$ $\delta>0$
and a map $f:(0,R)\times(0,\infty)\to \NN$
 there exists 
$\epsilon=\epsilon(n, R, \bar{R}, D, \delta, f)>0$ 
such that if $\bar{B}_p(D) \subset M^n$ is a closed ball in a complete
Riemannian manifold with the $f$ ball packing property
 such that
$M^n$ is locally $(\epsilon,R)$-almost isotropic off 
an $(\epsilon,R)$-almost unseen set of the form
described in (\ref{Wballs})
then 
\be \label{MainThmGH}
d_{GH}(B_p(D) \subset M^n, B_y(D)\subset Y)<\delta
\ee
where $Y$ is a complete length space with a subset
$W_Y=\{y_1, y_2,...\}$, such that $d_Y(y_j,y_k)\ge 2\bar{R}$ and such that
if $Y'$ is a connected component of $Y\setminus W_Y$
then the closure $Cl(Y')$ is isometric to an
$n$ dimensional space form with injectivity radius $>R$.

Furthermore a connected region of $M\setminus W_M$, $M'$, is uniformly 
$(\epsilon+\delta)$-almost isotropic everywhere,
\be \label{MainThmGH2}
d_{GH}(B_p(D)\cap Cl(M'), B_y(D)\cap Cl(Y'))<\delta.
\ee    
and
\be \label{MainThmexp}
|F_q(\theta, s,t)-F_K(\theta,s,t)|<\delta \qquad \forall \, q\in B_p(D)\cap M',   
\ee
where $F_K$ is the isotropy function of the space form $Cl(Y')$
and $K$ depends on $Y'$.
\end{theorem}

The basic idea here is that the space $Y$ is a space created by joining
together space forms at single points.  It is not a smooth manifold because
at the points where the space forms are joined there is no local chart.
Interestingly the space forms that are used to create $Y$ do not need
to have the same curvature:
$Y$ could be a sphere joined to a plane at a point.  It is possible that more than
two space forms are joined at a single point and that the space forms can be
joined at multiple points.  However only finitely
many may meet at any given, as can be seen since the limit space must be locally compact.
It is also possible that $Y$ is a hyperbolic space joined to
countably many spheres at countably many points.  It is possible that $Y$ could
be a single space form with two points that are set equal to one another
like a gateway.

Cosmologically, one can think of the points as black holes or gates or some unknown phenomenon
and the different space forms as being the universe as seen on various sides of
these points.

An example of a sequence of $M_i$ converging to $Y$ with a single bad point $W_Y$ 
is a sequence of rescalings of the Schwarzschild metric which converges
to a pair of planes joined at a point.  
It is easy to see how to extend this
example using cutting and paste techniques to give examples where
$Y$ is any space which is of the form described in Theorem~\ref{MainThm}.
See Lemma~\ref{schcut}, Corollary~\ref{schcutcor} and Example~\ref{anyY}.

Note that the uniform almost isotropy achieved in (\ref{ricciexp}) and (\ref{MainThmexp})
is not a consequence of the Gromov-Hausdorff closeness (\ref{ricciGH}) and (\ref{MainThmGH}).
These equations provide significant angular information about the Riemannian manifold
while  Gromov Hausdorff closeness can only be used to estimate distances.
It is quite possible that two manifolds be very close in the Gromov Hausdorff sense
and yet have very different formulas for the length of the third side of a triangle.
Consider the surface of a smooth ball versus the surface of a golf ball and the wild
behavior of geodesics on the latter.  
The proof (\ref{ricciexp}) and (\ref{MainThmexp}) involves an extension of Grove Petersen's 
Arzela Ascoli theorem \cite{GrPet}
and makes strong use of the almost isotropy condition.

%Cosmologically these two equationsich  are stating that the distances between points on geodesics
%which avoid strong lensing events can be estimated using the angle between the geodesics
%at their starting points and applying the standard formulas from the appropriate space form.
%Geodesics which pass through strong lensing events are not estimated in this manner.

%question were Such models were shown 
%by Dyer and Roeder to remain close to the Friedmann model as time
%evolves.  

%However it  also allows for
%significantly more irregular spacelike Riemannian manifolds and evolutions of
%spacetime should be investigated with these more general possibilities.
%For example, one might consider in what sense solutions of the Einstein
%equation vary continuously or not if the initial spacelike cauchy surfaces
%vary in the Gromov-Hausdorff topology.

We now provide a quick survey of the contents of this paper pointing out
key results which may be useful to mathematicians who study length spaces
and non-Euclidean geometry.  We also provide the definition of a space
we call and an {\em Exponential Length space} and relate it to the
above project. 

Section~\ref{sectlimits} has a review of Gromov Hausdorff theory.
In it we make the usual conversion of
Theorem~\ref{MainThm} into a theorem regarding
limits $Y$ of sequences of Riemannian manifolds $M_i$ 
which are locally $(\epsilon_i,R)$-almost isotropic
off $(\epsilon_i, R)$-almost unseen sets $W_{M_i}$
where $\epsilon_i \to 0$ [Theorems~\ref{MainThm2}].  
Ordinarily, even with an assumption of $Ricci \ge (n-1)H$, such a limit
space $Y$ is a complete length space with no well defined exponential
map.  Gromov proved that between every pair of points in the limit space
there is at least one length minimizing
curve which achieves the distance between the points \cite{Gr}.  However, it is
quite possible to have two length minimizing curves which overlap for
some time and then diverge.  This makes it very difficult to control what
happens to the $(\epsilon_i, R)$-almost isotropy in the limit process.
To prepare for this difficulty, we  extend the Grove-Petersen
Arzela-Ascoli Theorem \cite{GrPet} to a theorem concerning the limits
of almost equicontinuous functions on converging spaces
[Definition~\ref{almequicont} and Theorem~\ref{arzasc}].

Section~\ref{sectunseen} has examples of Riemannian manifolds
which are almost isotropic off almost unseen sets and
contains a couple of technical lemmas regarding such Riemannian manifolds
which are used in subsequent sections.

Section~\ref{beforeonto} studies the limiting behavior of the
exponential maps of the sequences $M_i$ proving Theorem~\ref{allbutonto}.
These limit exponential maps are not defined
on  a set $W_Y \subset Y$ [Definition~\ref{defWY}] but are homeomorphisms
onto their images and describe length minimizing curves.  
Note that in this section no assumption is
made that the almost unseen sets $W_{M_i}$ need to be contained in unions of
uniformly disjoint balls as in (\ref{Wballs}).

In Section~\ref{sectonto} we add this last condition (e.g.(\ref{Wballs}) with
$\epsilon=\epsilon_i$) on the $W_{M_i}$  and use it to prove that $W_Y$ is discrete
[Lemma~\ref{discrete}].
We then show that the exponential maps constructed on $Y$ are
locally surjective onto balls in $Y$ as long as they avoid $W_Y$.

At this point in the paper, enough properties of the limit space $Y$ will
have been proven to proceed and we will no longer need to refer to the sequence of
$M_i$.  Thus the remaining sections will be written about complete length spaces
that share these properties.  We make the following definition.

\begin{defn} \label{explsp} {\em
A complete length space $Y$ is called an {\em exponential length
space} off a set of points $W_Y$, if there exists a continuous function
$R_y>0$ and an {\em exponential dimension} $n\in \NN$ such that
for all  $y\in Y\setminus W_Y$, there is a continuous 1:1 function,
$exp_y:B_0(R_y) \to B_y(R_y)$, and 
\be \label{ryonto}
\exists r_y=d_Y(y, W_Y)\in (0,R_y] 
\textrm{ such that } exp_y:B_0(r_y) \to B_y(r_y)
\textrm{ is onto.}
\ee
Furthermore for fixed $v\in S^{n-1}$,
$exp_y(tv)$ is a length minimizing
curve for $t\in [0,R_y)$.

If $inf_y(R_y)= R>0$ exists we call $R$ the {\em exponential radius}.  }
\end{defn}

It should be noted that there is no assumption that the
exponential functions $exp_x$ vary continuously in 
the variable $x$ and that this allows us to avoid the issues
involved in defining a tangent bundle.
Note that the exponential radius plays a role similar to the injectivity
radius of a Riemannian manifold.  Zhongmin Shen has informed the
author that complete Finsler spaces are also
exponential length spaces  and in fact the 
exponential map is a $C^1$ diffeomorphism in that case.  
Example~\ref{sphpl} is an exponential length space off
a single point. 

In Section~\ref{sectexpl} we prove a few lemmas concerning 
exponential length spaces $Y$ off discrete subsets $W_Y$
with positive exponential radius.  
In Lemma~\ref{ontoconn} we prove that if $ Y'$ is a connected component
of $Y\setminus W_Y$, then for any $y\in Y'$
the exponential map, $exp_y$, is a homeomorphism 
from $B_0(R)$ onto $B_y(R) \cap Cl(Y')$.  The proof of this theorem
involves the Invariance of Domain Theorem, a strong topological result stating
that a subset of $\RR^n$ which is homeomorphic to a ball in $\RR^n$ is an open
subset of $\RR^n$ (c.f. \cite{EilSt}).

In Section~\ref{sectliso} we add in the condition of local isotropy which was
proven to hold on the limit spaces $Y$ of the $M_i$ in Theorem~\ref{allbutonto}.
Once again we make a definition to describe such spaces.

\begin{defn}\label{explsp2} {\em
An exponential length space $Y$  is {\em locally isotropic} off $W_Y$
if it has exponential radius greater than $R>0$ and for all 
$x\in Y \setminus W_Y$ 
there is a function $F_x: [0,\pi]\times [0,R]\times[0,R]$ 
which is continuous and satisfies 
\begin{eqnarray} 
a) & F_x(\theta_1, a, b)\le F_x(\theta_2, a,b) \quad 
\forall \theta_1<\theta_2 \label{abc2a}\\
b) & F_x(\pi, t,t) \ge t \label{abc2b}\\
c) & F_x(0,0,R)=R \label{abc2c}
\end{eqnarray}
such that
\be \label{aboiso}
d_Y(exp_x(tv),exp_x(sw))=F_x(d_S(v,w), t,s) \qquad \forall v,w\in S^{n-1},
\,\,\forall s,t\in [0,R).
\ee

We will call $R$ the {\em isotropy radius} and $F_x$ the isotropy function
about $x$.}
\end{defn}

Note that (\ref{abc2a})-(\ref{abc2c}) and 
(\ref{aboiso}) are the natural limits
of (\ref{abc1a})-(\ref{abc1c}) and (\ref{Fpdef})
in the definition [Defn~\ref{almiso}] of
a locally  $(\epsilon,R)$-almost isotropic Riemannian manifold as
$\epsilon$ is taken to $0$.

In Sections~\ref{sectliso} through~\ref{sectbirk}
we prove the following theorem and its corollary.

\begin{theorem} \label{isoexplthm}
If $Y$ is a locally isotropic
exponential length space off of a discrete set $W_Y$ 
then the closure of any connected component $Y'$ of $Y\setminus W_Y$
is a space form and its exponential structure
matches that of the space form.  In particular
$F_y(\theta,s,t)=F_K(\theta,s,t)$ for all $y\in Y'$
where $K$ is the sectional curvature of that space form.  
\end{theorem}

\begin{corollary}\label{isoexplcor}
If $Y$ is an exponential length space which is locally isotropic
then it is a complete Riemannian manifold with constant sectional
curvature. %fixed typo Mar 2004
\end{corollary}

Theorem~\ref{isoexplthm}
implies Theorem~\ref{MainThm} when combined with the
work in the earlier sections.  
The proof of Theorem~\ref{isoexplthm} is broken into sections
which focus on different properties of exponential length spaces.
We hope that they will prove interesting to mathematicians who
study length spaces and non-Euclidean geometry.

In Section~\ref{sectliso} we first prove that a locally isotropic exponential
length space is  uniformly locally minimizing [Definition~\ref{defulocalmin}, Lemma~\ref{localmin}].
That is there is a uniform distance such that pairs of points which are less than that
distance apart have unique length minimizing curves running between them. 
We then construct local isometries
between balls centered in $Y'$ [Lemma~\ref{soisom} and Corollary~\ref{Fconst2}]
and then show that $Cl(Y')$ 
is isometric to a locally isotropic exponential length space 
everywhere whose isotropy functions do not depend upon the
base point [Lemma~\ref{badsexp}].  In some weak sense we now know have that $Cl(Y')$
has constant sectional curvature but we do not yet have smoothness.

In Section~\ref{sectexpcurve} we prove that uniformly locally minimizing
exponential length spaces are extended exponential length spaces [Definition~\ref{extexplsp}].
That is, the exponential maps on these spaces are not just defined on balls but
can be extended to maps on all of $\RR^n$ with good properties.
There is no assumption of local isotropy made in this section.
To extend the exponential map continuously to
all of $\RR^n$  we introduce the concept of {\em exponential curves}
[Definition~\ref{defexpcurve}, Lemma~\ref{extend}
and Lemma~\ref{extcont}].  
We show  that in such spaces,
all length minimizing curves are exponential curves
[Lemma~\ref{minexp}] and all exponential curves are locally length minimizing.  

In Section~\ref{sectext} we study the concept of a {\em conjugate point} 
on an {\em extended exponential length space}.
[Definition~\ref{extexplsp} and Definition~\ref{conjpt}].  
Through a series of lemmas, we show that, in such a space,
if a ball has no conjugate points
then it is mapped by the exponential map as a local homeomorphism
[Theorem~\ref{conjhomeom}].  
This extends the traditional theorem in Riemannian geometry which 
obtains a local diffeomorphism using the Inverse Function Theorem
in regions without conjugate points.
We apply the Invariance of Domain Theorem here to replace the
Inverse function theorem because we do not yet have smoothness.

In Section~\ref{secttube} we reintroduce the hypothesis of local isotropy.
We prove that we can extend isometries between subsets to isometries of balls [Lemma~\ref{ballext}]
and use them to prove that the distances between points on closely
located pairs of exponential curves depends only on the angle between them [Lemma~\ref{ctogether}].
We then show that the distance to conjugate points is a constant on $Y$
[Lemma~\ref{allconj}].

In Section~\ref{simpconn} we add the condition that $Y$ is simply connected.
We then prove $Y$ is homeomorphic to $\RR^n$, $\SSS^n$ 
or $\HH^n$ via the exponential map [Theorem~\ref{homeomdone}].  In Lemma~\ref{triangle} we extend
Lemma~\ref{solveF} to triangles of all sizes.  Then we construct global isometries
in Lemmas~\ref{globisom} and~\ref{isomext}.

In Section~\ref{sectbirk} we complete the proof
of Theorem~\ref{isoexplthm}.  The key ingredient is Birkhoff's Theorem
which states that if a space $X$ has locally unique length minimizing curves
and any isometry on subsets of $X$ extends to a global isometry
then the space must be $\SSS^n$, $\HH^n$ or $\EE^n$ \cite{Bi}.
We then prove Lemma~\ref{match} matching the
exponential structure of the length space to that of the space form
and complete the proof of Theorem~\ref{MainThm2}.

Finally in Section~\ref{sectricci} we prove Theorem~\ref{RicciThm}
by demonstrating that when a lower bound on Ricci curvature is assumed one
cannot have a limit space which contains pairs of space forms joined
at a point.  The proof consists of a careful measurement of the
volumes of balls using the Bishop-Gromov Volume Comparison Theorem.
Recall that such a comparison holds on the limit space $Y$ by
Colding's Volume Convergence Theorem \cite{Co}.

The author would like to thank Gregory Galloway  for
being a guide to general relativity,  Jeffrey Weeks for his inspiring 
presentation
on the Friedmann model at Lehigh University and P.J.E. Peebles for his excellent webpage and textbook
\cite{Peeb}.  She would especially like to thank Richard Schoen for suggesting a
stronger conclusion should hold for Theorem~\ref{RicciThm} than Theorem~\ref{MainThm}.  
She would like to thank Al Vasquez for the reference
to the Invariance of Domain Theorem and mathscinet for producing Birkhoff's
paper when searching with the words  ``spherical'', ``hyperbolic'' and ``characterized''.
She would like to thank Professors Cornish, Uzan, Souradeep and Labini for
references and information about cosmology and Ivan Blank and Penelope Smith
for discussions regarding the concept of almost unseen sets.

%Weeks, Cornish,  Uzan, Souradeep, Labini
%for references and info about cosmology and the anisotropy of the
%universe.

\sect{Limits} \label{sectlimits}   

In this section we review the definition
of Gromov-Hausdorff convergence, Gromov's compactness theorems \cite{Gr}
and Grove-Petersen's Arzela-Ascoli Theorem \cite{GrPet} and some 
other key concepts.  
We extend this Arzela Ascoli Theorem to functions which are only
almost equicontinuous [Definition~\ref{almequicont} and Theorem~\ref{arzasc}].
Finally we reduce Theorem~\ref{MainThm} to the
following theorem.

\begin{theorem} \label{MainThm2}
Given  $n \in \NN$, $\bar{R}, R>0$, $\delta>0$
and a map $f:(0,R)\times(0,\infty)\to \NN$, $\epsilon_i \to 0$,
if $M_i^n$ are complete
Riemannian manifolds with the $f$ ball packing property
 such that
$M_i^n$ is locally $(\epsilon_i,R)$-almost isotropic off 
an $(\epsilon_i,R)$-almost unseen set $W_{M_i}$ of the form
described in (\ref{Wballs}) with $\epsilon=\epsilon_i$ respectively
then a subsequence of the $M_i$ converges to a complete length space
$Y$ with a countable collection of points
$W_Y=\{y_1,y_2,...\}$, such that $d_Y(y_j,y_k)\ge 2R$ and such that
if $Y'$ is a connected component of $Y\setminus W_Y$
then $Cl(Y')$ is isometric to an
$n$ dimensional space form.

Furthermore if
$M_i'$ are connected regions of $M_i\setminus W_{M_i}$
containing points $q_i$ converging to $x\in Y$, 
then $GH\lim_{i\to\infty} (Cl(M'_i),q_i)=(Cl(Y'),x)$
For fixed $s,t \in [0,R)$, $\theta\in [0,\pi]$,
the almost isotropy functions converge,
\be \label{MainThm2a}
\lim_{i\to\infty}F_{q_i}(\theta, s,t)=F_K(\theta,s,t)
\ee
where $K$ is the sectional curvature of $Cl(Y')$.
\end{theorem}

\begin{remark}\label{soweirmk}
Note that in \cite{SoWei}, it is proven that if the $M_i$ are
compact and simply connected, then the limit space $Y$ is its own universal cover.
So in that case $Y'$ is not just a space form but it must be
$\SSS^2$, as the sphere is the only compact simply connected space form.  
Thus compact simply connected $M$ in Theorems~\ref{MainThm} and~\ref{RicciThm}
are close to spaces where $Y'$ and $Y$ are spheres respectively.
\end{remark}

Those who are experts in
this theory will immediately see why this theorem implies
Theorem~\ref{MainThm} and can skip ahead to the next section.
The discussion of equicontinuity [Defn~\ref{almequicont} and Thm~\ref{arzasc}]
can be referred to as needed later on.

We now provide the necessary background for following the remainder of this paper.
We recommend \cite{BBI} as a reference for non-experts.  

\begin{defn} \label{ms}
A {\em metric space} is a set of points, $X$, and a distance function $d:X\times X \to [0,\infty)$ 
such that $d(x,y)=0$ iff $x=y$, $d(x,y)=d(y,x)$ and $d(x,y)+d(y,z)=d(x,z)$ for all $x,y,z \in X$.
\end{defn}

\begin{defn} \label{cls} 
A {\em complete length space} is a metric space such that between every pair
of points there is a length minimizing (rectifiable) curve joining them whose length is
the distance between the points.
\end{defn}

%\begin{defn} \label{Haus} {\em [Hausdorff]
%If $X_1$ and $X_2$ are metric spaces with isometric embeddings $f_i:X_i \to Z$ then
%their {\em Hausdorff distance} is defined
%\be
%d_Z(f_1(X_1), f_2(X_2))= \inf \{ \epsilon: T_\epsilon(f_1(X_1))\supset f_2(X_2) \textrm{ and } T_\epsilon(f_2(X_2))\supset f_1(X_1)\}.
%\ee }
%\end{defn}
%
%The idea is that $X_1$ and $X_2$ are close to each other when they are embedded in $Z$ but not
%that they have similar shapes.  For example if $f_1(X_1)$ is the circle $\{(x,y,z):x^2+y^2=10, z=0\}$
%and $f_r(X_r)$ is the torus $\{(x,y,z): (\sqrt{x^2+y^2}-10)^2 +z^2=r^2\}$ where $r<1$ then 
%\be
%d_Z(f_1(X_1), f_2(X_r))=r
%\ee 
%and as $r$ converges to $0$, $f_r(X_r)$ converges to $f_1(X_1)$ in the Hausdorff sense.
%
%To compare metric spaces which are not embedded in a common space, Gromov made the following defintion.
%
%\begin{defn} \label{Ghclose}  {\em [Gromov]
%If $X_1$ and $X_2$ are metric spaces then their {\em Gromov -Hausdorff distance} is defined
%\be
%d_{GH}(X_1, X_2) =\inf \{ d_Z(f_1(X_1), f_2(X_2)): \textrm{ all metric spaces } Z \textrm{and all isom embeddings } f_i: X_i \to Z \}
%\ee  }
%\end{defn}

Gromov's definition of the Gromov-Hausdorff distance between two spaces involves infima
of the Hausdorff distances between all possible embeddings of these spaces.  We will not
be using this definition but rather a very useful property of the Gromov-Hausdorff
distance which relates the concept to maps between the two spaces.

\begin{defn} \label{GH1} {\em
A function $\phi:X \to Y$ is an {\em $\epsilon$-almost isometry} if it is 
\be
\epsilon \textrm{-almost
distance preserving: } |d_Y(\phi(x_1), \phi(x_2))-d_X(x_1, x_2)|<\epsilon  
\ee
and
\be
\epsilon \textrm{-almost onto } T_\epsilon(\phi(X))\supset Y.
\ee
Note that $\phi$ need not be continuous.}
\end{defn}

Recall that two spaces are isometric when there is a map, called an isometry, between them which is 1:1, onto and distance
preserving.  Actually the fact any map that is distance preserving is 1:1.
So an almost isometry is an approximation of the concept of an isometry.
The fact that the almost isometry is neither continuous nor onto not 1:1 allows the two spaces to be shaped
quite differently.  For example $X$ could be a circle and $Y$ could be a thin torus and we would
get an almost isometry from $X$ to $Y$ by embedding $X$ in $Y$ and an almost isometry from
$Y$ to $X$ by mapping rings to single points.  Notice that both the topology and the dimensions
of $X$ and $Y$ are quite different.  

The following lemma [c.f. \cite{BBI}], will be used in place of a definition both when
proving and when applying Theorem~\ref{MainThm}.
 
\begin{lemma} \label{GH2} 
Suppose $X$ and $Y$ are metric spaces then
\be
\textrm{ if }d_{GH}(X,Y) < \epsilon \textrm{ then there is a } (2\epsilon) 
\textrm{-almost isometry from } X \textrm{ to } Y,
\ee 
and
\be
\textrm{ if there is an } \epsilon 
\textrm{-almost isometry from } X \textrm{ to } Y
\textrm{ then }d_{GH}(X,Y) < 2\epsilon.
\ee
\end{lemma}

If it were not for the annoyance of the change of $\epsilon$ to $2\epsilon$
the existence of an almost isometry between 2 spaces would make a wonderful 
definition of Gromov-Hausdorff distance. 

Cosmologically,
Theorem~\ref{RicciThm} then says that we have an almost isometry between an almost isotropic manifold, $M$
and a space form $Y$.  This means that distances between points can be estimated using an almost isometry to the space
form.  The actual almost isometry is not produced in this paper which makes this difficult to apply
cosmologically.  Nevertheless there are implications.  Space forms have lots of isometries.  In fact
all balls whose radius is less than the injectivity radius are isometric to each other.  
Using the almost isometry between $M$ and
$Y$ we get the fact that all balls of this size are almost isometric in $M$, that space looks pretty
much the same from point to point.  This is a stronger fact than (\ref{ricciexp}) because some
balls may contain components of $W$ where there is strong gravitational lensing and now one
can estimate the distances between stars which cannot see each other without passing through $W$.
In fact, one can use the region of space near earth as a sample ball (which does not contain
any strong gravitational lensing) and then know that distant regions (even those containing
strong gravitational lensing) are almost isometric.  The actual black hole would have to
be on the scale of the error, $\epsilon$, in the almost isometry but on the cosmological scale
things could be well understood.  It is often assumed that regions around black holes look
just like Euclidean space, here we have proven that they must be close to a space form but
not necessarily in a smooth way.  
Theorem~\ref{MainThm} essentially has the same result where we compare $M'$ the connected component
of $M$ (the part of space which can be reached without passing through $W$) to a space form $Y'$.

Once one has an understanding of Gromov-Hausdorff distance, one can define 
the convergence of metric spaces.  That is metric spaces $X_i$ converge to a metric space
$Y$ iff $d_{GH}(X_i, Y)$ converges to $0$.  This definition is too
restricted for applications with unbounded limit spaces so Gromov defined the following 
pointed Gromov Hausdorff convergence.

\begin{defn} \label{ptGHconv} 
{\em [Gromov]
If each $x_i$ is in a complete metric space $X_i$, we say $(X_i,x_i)$ converges to $(X_0, x_0)$ in the {\em pointed
Gromov Hausdorff sense} if for all $D>0$  the closed balls $B(x_i,R)\subset X_i$ converge
in the Gromov Hausdorff sense to $B(x_0,D)\subset X_0$.}
\end{defn}

He then proved the Gromov Compactness Theorem:

\begin{theorem} \label{GrCompactness} {\em [Gromov]}
If $X_i$ are complete length spaces that satisfy a uniform $f$ 
ball packing condition, then 
for any $x_i \in X_i$ a subsequence of $(X_i,x_i)$ converges to a complete length space $(Y,y)$
in the pointed Gromov-Hausdorff sense.  Conversely, if $(X_i,x_i)$ converge to a complete
length space $(Y,y)$ then they satisfy a uniform $f$ ball packing condition.
\end{theorem}

As a consequence any sequence of complete Riemannian manifolds 
with a uniform lower bound
on Ricci curvature converges to a complete length space.  Note, however, that in general the
limit space will not be a manifold.  For example a sequence of hyperboloids can converge to a
cone and a sequence of paraboloids to a half line.  

Note also that the closeness in Gromov's
compactness theorem is on compact regions, not on the whole manifold at once.  This is why
there are balls $B_p(D)$ mentioned in the Theorem's~\ref{MainThm} and~\ref{RicciThm}.
This condition is necessary in the statement of these theorems as can be seen
in the following example.

\begin{example} \label{ballchange}
{\em
Suppose $M_s^2$ are warped product manifolds with the metrics,
\be
g= dt^2 + sinh^2(K_s(t) t) g_0,
\ee
where $K_s(t)$ is increasing and smooth such that
$K_s(t)=1$ on $[0,1]$ and $K_s(t)= 1+ (Ln(t))^{1/3}(1/s) $ on $[2, e^{27s^3}]$
and $K_s(t)= 4+1/(se^{27s^3}+s)$ on $[e^{27s^3}+1, \infty) $.
Note that this can be done smoothly keeping $K'_s(t)$ and $K''_s(t)$ both less than 
$5/s$ on $[0,2]$ and 
both less than $5/ (se^{27s^3}+s)$ on $[e^{27s^3},e^{27s^3}+1]$.

On a 2 dimensional warped product, the sectional curvature
for $q\in\partial B_p(t)$ is 
\begin{eqnarray*}
Sect_q &= &-\frac{ \frac{\partial^2}{\partial t^2} (sinh(K_s(t) t) }{(sinh(K_s(t) t)}\\
& = &-\frac{ \frac{\partial}{\partial t} (cosh(K_s(t) t) (K_s'(t)t+K_s(t)) }{(sinh(K_s(t) t)} \\
&= & -\frac{  (sinh(K_s(t) t) ( K_s'(t)t+K_s(t) )^2 }{(sinh(K_s(t) t)} \\
  & \,\,& -\frac{  (cosh(K_s(t) t) ( K_s''(t)t+ K'_s(t)+K_s'(t) ) }{(sinh(K_s(t) t)} \\
&=& - (K_s(t)-K_s'(t)t)^2
  - ( K_s''(t)t+ K_s'(t)+K_s'(t) ) (coth(K_s(t) t)). \\
&=& - (K_s(t))^2 + 2K_s(t)K_s'(t)t - (K_s'(t)t)^2
  - ( K_s''(t)t+ 2K_s'(t) ) (coth(K_s(t) t)). \\
\end{eqnarray*}

Suppose we fix a number $R>0$.  For any $p\in M_s^2$, let 
$r=d(p, p_0)$ then  $q\in B_p(R) \subset Ann_{p_0}(r-R, r+R)$,
so $t=d(q, p_0)\in (r-R, r+R)$.  Then the sectional curvature at $q$ is close to a constant
$-(K_s(r))^2$ for sufficiently large $s$ as follows:
\begin{eqnarray*}
|Sect_q+K_s(r)^2|&\le&|K_s(t)^2-K_s(r)^2|+| 2K_s(t)K_s'(t)t| +| (K_s'(t)t|^2
  + | K_s''(t)t| + |2K_s'(t)| \\
&\le&2 \max_{a\in [r-R, r+R]}( 2K_s(a)K_s'(a) R ) \\
& \,\, & + | 2K_s(t)K_s'(t)t| +| (K_s'(t)t|^2
  + | K_s''(t)t| + |2K_s'(t)|,
\end{eqnarray*}
which can be shown to be small by examining the following sets of cases:

First, we have
\begin{eqnarray*}
 2K_s(a)K_s'(a) R  & = &0 \textrm{ for } a\in[0,1], \\
 2K_s(a)K_s'(a) R&\le &  2(1+(Ln(2))^{1/3}(1/s))(5/s)R \textrm{ for } a\in[1,2], \\
 2K_s(a)K_s'(a) R&\le &  |2(1+(Ln(a))^{1/3}/s)( (1/3)(Ln(a))^{-2/3}/(as) )R \\
&\le&  |2(1+(Ln(2))^{1/3}/s)( (1/3)(Ln(2))^{-2/3}/(2s) )R
\textrm{ for } a\in[2, e^{27s^3}] \\   
 2K_s(a)K_s'(a) R &\le &  2(5)(5/(se^{27s^3}+s))(R)  \textrm{ for } a\in 
[e^{27s^3}, e^{27s^3}+1], \\
2K_s(a)K_s'(a) R&\le & 0 \textrm{ for } a\in [e^{27s^3}+1, \infty).
\end{eqnarray*}
Then we bound 
\begin{eqnarray*}
| 2K_s(t)K_s'(t)t| &+& | (K_s'(t)t|^2
  + | K_s''(t)t| + |2K_s'(t)| \le   \\
\,\,&\le&  0   \textrm{ for } t\in[0,1], \\
| 2K_s(t)K_s'(t)t| &+& | (K_s'(t)t|^2
  + | K_s''(t)t| + |2K_s'(t)| \le   \\
&\le&  2(1+(Ln(2))^{1/3}(1/s))(5/s)2
      + 100/s^2+10/s + 10/s \textrm{ for } t\in [1,2], \\
| 2K_s(t)K_s'(t)t| &+& | (K_s'(t)t|^2
  + | K_s''(t)t| + |2K_s'(t)|   \le  \\
&\le&   |2(1+(Ln(t))^{1/3}/s)( (1/3)(Ln(t))^{-2/3}/(ts) )t 
       +( (1/3)(Ln(t))^{-2/3}/(s) )^2    \\
&\,\,&      +| -2(Ln(t))^{-5/3}/(9ts) - (Ln(t))^{-2/3}/(3ts) |
   +|2(Ln(t))^{-2/3}/(3ts) )|    \\
&\le &  |2(1+(Ln(2))^{1/3}/s)( (1/3)(Ln(2))^{-2/3}/(s) ) +
   ( (1/3)(Ln(2))^{-2/3}/(s) )^2  \\
& \,\,&\qquad      +| (-2)
(Ln(2))^{-5/3}/(18s) - (Ln(2))^{-2/3}/(6s) | \\
&\,\,   &\qquad
+|2( (1/3)(Ln(2))^{-2/3}/(2s) )| 
   \textrm{ for } t\in[2, e^{27s^3}],   \\
| 2K_s(t)K_s'(t)t| &+& | (K_s'(t)t|^2
  + | K_s''(t)t| + |2K_s'(t)| \le   \\
&\le &  2(5)(5/(se^{27s^3}+s))(e^{27s^3}+1)+ 
(5/(se^{27s^3}+s))^2(e^{27s^3}+1)^2+ \\
&\,\,&\,\, + (5/(se^{27s^3}+s))(e^{27s^3}+1)
+ 2(5/(se^{27s^3}+s))   
 \\
 &\,\,&\,\,
+ 100/s^2+10/s + 10/s   \textrm{   for } t\in[e^{27s^3}, e^{27s^3}+1],  \\
| 2K_s(t)K_s'(t)t| &+& | (K_s'(t)t|^2
  + | K_s''(t)t| + |2K_s'(t)| \le   \\
&\le & 0  \textrm{ for } t\in [e^{27s^3}+1, \infty).
\end{eqnarray*}

That is, for all $R>0$ and $\epsilon>0$, there exists $s$ sufficiently large
that $|sect_q-K_s(d(p,p_0))^2|<\epsilon$ for all $q, p \in M_s^n$ such that $d(p,q)<R$.
Since the distance between geodesics emanating from $p$ can be estimated from
above and below by integrating the curvature, this implies that there exists
$s$ sufficiently large depending only on $R$ and $\epsilon'$ such that
\be
|d_{M_s}(exp_p(tv), exp_p(sw))-F_{K_s(d(p,p_0))^2}(d_S(v,w),t,s)|<\epsilon'
\ee
for all $t,s <R$, for all $p\in M_s^n$. 

Since the $M_s$ also have curvature uniformly bounded below by $-25$, they satisfy the 
conditions of Theorem~\ref{RicciThm} and so balls of a fixed radius
$D$ approach a space form.  However different balls in $M_s$ will 
approach different space forms.  In particular, the ball near the
center of $M_s$ is a space of constant curvature $1$ while a ball
far away from the center will have constant curvature greater than $16$. }
\end{example}

Now it is common to refer to the concept of points $p_i$ in the $X_i$ converging
to a point $z$ in the limit $Y$.  This is made rigorous if one uses the 
almost isometries $\phi_i:X_i \to Y$
from Lemma~\ref{GH2}.  We first choose the isometries $\phi_i$ to fix a particular convergence onto
the limit space.  For example when a hyperboloid converges to a cone the $\phi_i$ can be rotated many different
ways.  We need to fix the $\phi_i$ to discuss particular points.
Then we say $p_i$ converge to $z$ if $\phi_i(p_i)$ converge to $z$ as points in $Z$.
Note that given any sequence of $p_i \in B(x_i,D)$ we know a subsequence converges because $B(y,D)$ is compact.

We can now prove that Theorem~\ref{MainThm2} implies Theorem~\ref{MainThm}.

\vspace{.5cm}

\noindent{\bf Proof:}
Suppose on the contrary that Theorem~\ref{MainThm} is false
for some $n \in \NN$, $\bar{R}, R>0$, $D>0$ $\delta>0$
and a map $f:(0,R)\times(0,\infty)\to \NN$.  So 
there is a sequence of  $\epsilon_i$ converging to $0$
and a sequence of manifolds $M_i$ which 
satisfy the $f$ ball packing property and are
locally $(\epsilon_i,R)$-almost isotropic off 
$(\epsilon_i,R)$-almost unseen sets, $W_i$ of the form
described in (\ref{Wballs}) with $\epsilon=\epsilon_i$
such that for all $i$,
(\ref{MainThmGH}) and (\ref{MainThmexp}) don't hold for
any complete length space $Y$ of the form described in the
theorem.  However Theorem~\ref{MainThm2} states that they must
converge in the pointed Gromov-Hausdorff sense to exactly such a space
$Y$, which means that for $i$ sufficiently large, depending on
$\delta$, we do in fact have a length space $Y$ satisfying
(\ref{MainThmGH}) and (\ref{MainThmexp}) which contradicts
the above.
\qed

\begin{example} \label{packreq}
The ball packing condition is necessary in Theorem~\ref{MainThm2} because
one can take a sequence of hyperbolic manifolds with constant curvature
$K$ and take $K$ to negative infinity.  Each space is actually isotropic
but the sequence does not have a converging subsequence approaching a 
length space $Y$.
\end{example}

It is an open question as to whether the ball packing condition is
necessary in Theorem~\ref{MainThm}.  To try to prove this theorem
without the ball packing condition would involve adapting the
Gromov Compactness Theorem to say something about sequences
of manifolds which don't have converging subsequences, a daunting task.

Clearly the first step towards proving Theorem~\ref{MainThm2}
will be to apply Gromov's Compactness Theorem to obtain a limit space $Y$.
However, to study the isotropy on $Y$ we will need an exponential
map.  We will construct such an exponential map by taking the limit
of a subsequence of the exponential maps defined on $M_i$.  

There is already
an extension of Arzela Ascoli Theorem to Gromov Hausdorff situations
by Grove and Petersen \cite{GrPet},
which states that if a sequence of continuous functions $f_i:X_i\to Y_i$
are equicontinuous and $(X_i,x_i) \to (X,x)$ and $(Y_i,f(x_i)) \to (Y,y)$ 
in the pointed Gromov Hausdorff sense
then a subsequence of the $f_i$ converge to a limit function
$f:X\to Y$.  This implies that curves which are parametrized by arclength converge
and that length minimizing curves converge to length minimizing curves but
it does not control the angular behavior of the exponential maps.

\begin{figure}[htbp]
\includegraphics[height=2.2in ]{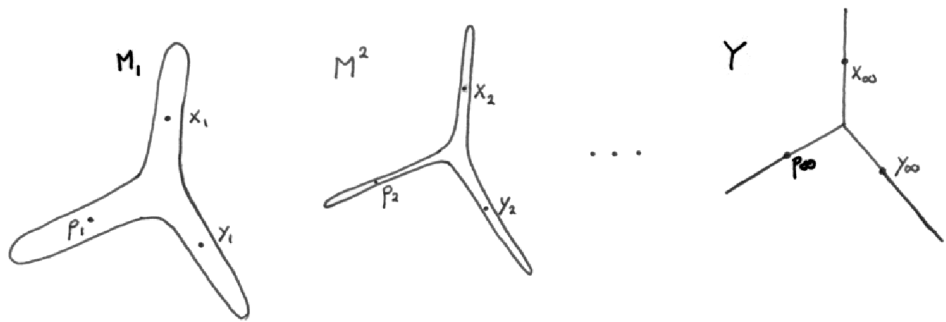}
\caption{ Pairs of geodesics running from $p_i$ to $x_i$ and from $p_i$ to $y_i$ converge
to a pair of minimizing curves in $Y$ running from $p_\infty$ to $x_\infty$ and from
$p_\infty$ to $y_\infty$.  This limit pair start as an identical curve and then diverge. } \label{figY}
\end{figure}

In general the exponential maps are not well controlled under Gromov-Hausdorff
convergence.  For example, take a length spcae $Y$ consisting of 3 line segments meeting
at a point.  Suppose we have a sequence of Riemannian surfaces $M_i$
which converge to a $Y$ shaped $Y$.  Note that the exponential maps 
must converge to functions which are no longer injective and that there 
are minimizing curves which diverge from one another after initially 
overlapping.  See Figure~\ref{figY}.

%\psdraft
%\psfull

In this paper we can use almost isotropy to  control the exponential maps
to some extent and this will allow us to create better exponential maps 
on our limit spaces.  To do so we need the following more general
theorem which will allow the maps $f_i$ not to be continuous.

\begin{defn} \label{almequicont} {\em
A sequence of functions between compact metric spaces,
$f_i:X_i\to Y_i$, is said to be {\em uniformly almost
equicontinuous} if there exists $\epsilon_i$ decreasing to $ 0$ such that
for all $\epsilon>0$ there exists $\delta_\epsilon>0$ such that
\be
d_{Y_i}(f_i(x_1), f_i(x_2))<\epsilon+\epsilon_i
 \textrm{ whenever } d_{X_i}(x_1, x_2)<\delta_\epsilon.
\ee  }
\end{defn} 

\begin{theorem} \label{arzasc}
If $f_i:X_i\to Y_i$, is uniformly almost equicontinuous
between
complete length spaces
 $(X_i,x_i)\to (X,x)$ and $(Y_i, f_i(x_i))\to (Y,y)$ 
converge in the Gromov Hausdorff
sense where $X$ and $Y$ are compact, then a subsequence
of the $f_i$ converge to a continuous limit function
$f:X\to Y$.
\end{theorem}

\Pf
As in Petersen, choose countable dense subsets $A_i=\{a_1^i,a_2^i...\}\subset X_i$,
such that $a_j^i\to a_j\in X$ where $A=\{a_1, a_2,...\}$ is dense in $X$.
Then subsequences of $f_i(a_j^i)$ converge using the pointed
convergence of the $Y_i$ and the precompactness of balls
in the $Y_i$.  So we can thus apply the standard diagonalization argument to
get a subsequence of the $f_i$ which converges on these countable
dense sets to some function $f:A \to Y$.

We need only show $f$ is continuous on $A$ and then we can extend it
to a continuous function on $X$.  For all $\epsilon>0$ take
$N$ sufficiently large that $\epsilon_N<\epsilon/2$ and
$\delta<\delta_{\epsilon/2}$ so
\be
d_{Y_i}(f_i(x_1), f_i(x_2))<\epsilon
 \textrm{ whenever } d_{X_i}(x_1, x_2)<\delta, \,\,\, i\ge N.
\ee
So now given $a_j,a_k\in A$ such that
$d_X(a_j,a_k)<\delta/2$, taking $i\ge N$ sufficiently large that
$d_{X_i}(a^i_j,a^i_k)<\delta$, so
\be
d_{Y_i}(f_i(a^i_k), f_i(a^i_j))<\epsilon.
\ee
Then taking $i\to \infty$ we get
\be
d_{Y}(f(a_k), f(a_j))<\epsilon.
\ee
\qed

%------------------------------
\sect{Almost Isotropy off Almost Unseen Sets} \label{sectunseen}

In this section we provide some examples of Riemannian manifolds which
are almost isotropic off almost unseen sets [Definition~\ref{unseen}].  We also prove two
technical lemmas regarding such Riemannian manifolds which will be needed
later.

\begin{lemma}\label{unseenball}
A ball of radius $r$ in a space form $N$ of constant sectional
curvature $K$ will be $(\epsilon, R)$-almost
unseen if 
\be \label{unseenball1}
r<F_K(\epsilon, \epsilon, \epsilon)
\ee
 where
$F_K$ is the almost isotropy function of $N$ and $R<$ the
injectivity radius of $N$.
\end{lemma}

\Pf
If $p \in M_r\setminus T_\epsilon(B_q(r))$ then by our choice of $R$
and the symmetry of space forms $S_p$ is exactly one ball.  
Now let $v$ be such that $exp_p(d(p,q)v)=q$ and let $w \notin S_p$.
Then 
\be
F_{K}(d_S(v,w),d(p,q), d(p,q))\ge r,
\ee
and the radius, $\theta$,  of $S_p$ must satisfy
\be
F_{K}(\theta,d(p,q), d(p,q))\ge r,
\ee
as well.   So  by (\ref{unseenball1}),
\be
F_{K}(\theta,\epsilon, \epsilon)<F_{K}(\theta,d(p,p_i), d(p,p_i))=r<F_{K}(\epsilon,\epsilon, \epsilon),
\ee
 which implies that $\theta <\epsilon$ by the monotonicity of $F$.
\qed

We can now apply this lemma to construct examples of spaces with 
almost unseen Schwarzschild necks
and interesting limit spaces $Y$.

\begin{lemma} \label{schcut}
For any fixed pair of  space forms $N_1^3$ and $N_2^3$ and points $p_i\in N_i$
and given any $r$ sufficiently small,
we can construct a Riemannian manifold $M_r$ which is isometric to
$N_i \setminus B_{p_i}(r)$ on two regions and has a
Schwarzschild neck of diameter less than $16(r+4r^2)$
so that taking $r$ to zero we get a sequence of Riemannian manifolds
satisfying the conditions of Theorem~\ref{MainThm2}
which converge to a length space $Y=N_1 \cup N_2\, / \,p_1\sim p_2$.
\end{lemma}

\Pf
Recall that on a ball of sufficiently small radius $r$
the metric in a space form $N_i$ may be described as a warped product
metric $dt^2+ f_K(t)^2 g_0$ where $t\in [0,r)$,
$g_0$ is the standard metric on the
sphere and $f_K$ is either $sinh(\sqrt{-K}t)$, $t$ or $sin(\sqrt{K}t)$.

Recall that the Schwarzschild metric on $\RR^3\setminus \{0\}$
can also be described with a warped product metric as
\be
g_{sch}=(1+\frac{m}{2t})^4 g_{Euch}=(1+\frac{m}{2t})^4 (dt^2+ t^2 g_0).
\ee
(c.f. \cite{SchYau} \cite{Br})

Setting $s=m^2/(4t)$, we see that we get an isometric inversion
so that the Schwarzschild solution is asymptotically flat in both
directions.  In fact we can describe how to glue
our Schwartzchild neck for $t\ge m/2$ into $N_1$ and
then repeat the process to glue in $s\ge m/2$ into $N_2$.

Given $r>0$ as above and such that $f_K(r)<2r$, let $m_r<r$ such that 
So for radii $t\ge m/2$ on $\RR^3$ define the metric
\be
g_r=h_r(t)(dt^2+f_r(t)^2 g_0)
\ee
where
$h_r(t)$ is a smooth function with values in
$[1, (1+ \frac{m_r}{2t})^4]$ that is $1$ for $t\ge r$ and $(1+\frac{m_r}{2t})^4$
for $t<r/2$ and  $f_r(t)$ is a smooth function with values between
$t$ and $f_K(t)$ which is $f_{K_1}(t)$ for $t\ge r$ and $t$ for $t\le r/2$.
Then $g_r$ and its corresponding metric to glue in $M_2$ define a metric for
a Schwarzschild neck (with $t\le r$ and $s\le r$) that can be glued smoothly
to $N_i\setminus B_{p_i}(r)$ to create $M_r$.  The diameter
of the neck is then bounded above by
\begin{eqnarray}
(\max h_r) & 2(r-m_r/2) + (\max h_r)(\max f_r)^2\pi = & \\
\, &= (1+r/(2r))^4 (2r-r)+(1+r/(2r))^4(2r)^2 \pi &\le 16(r+4r^2).
\end{eqnarray}

It is easy to verify that as $r$ decreases to $0$, $M_r$ converges in
the Gromov Hausdorff sense to $Y$ and thus by Gromov's Compactness Theorem,
they satisfy a uniform $f$ ball packing property.  Furthermore if we set
$W_{M_r}$ to be the neck, it is clearly contained in a ball of radius
$16(r+4r^2)$.  If we set $R$ less than the minimum of the two injectivity
radii of the $M_i$, then it is easy to see that $M_r$ is $(0, R)$-almost
isotropic off $W_{M_r}$.  

Lastly taking any $\epsilon>0$ and 
setting $r$ according to the isotropy functions of the $N_i$ as follows,
\be \label{schcutr} 
r < \min_{i=1,2} F_{K_i}(\epsilon,\epsilon, \epsilon),
\ee
 we can verify that $w_{M_r}$ is
$(\epsilon,R)$-almost unseen from each $N_i$ by using the fact that any
geodesic entering the neck from $N_i$ passes into $B_{p_i}(r)$
and then applying Lemma~\ref{unseenball}.

\qed

Note that if we take $N_i$ to be spheres we could even choose $f_r$ and $h_r$
so that $M_r$ had nonnegative scalar curvature.

\begin{coro} \label{schcutcor}
Given any countable collection of space forms $N_i^n$ and points $p_{i,j}\in N_i$,
such that $d_{M_i}(p_{i,j}, p_{i,k})=d_{i,j,k}> 2\bar{R}$
and a bijective map without fixed points $P: \{ p_{i,j}\} \to \{ p_{i,j}\}$, 
we can construct a sequence of smooth Riemannian manifolds $M_r$ 
(locally as above) satisfying
the conditions of Theorem~\ref{MainThm2} which converge to the space 
\be
Y=\bigcup N_i \,\, / \,p_{i,k}~P(p_{i,k})
\ee
This includes the possibility of a single $N_i$ with an even number of points.
\end{coro}

\begin{example} \label{anyY}
In order to construct limit spaces $Y$ with many space forms meeting at a single point
we don't use the Schwartzchild metric as a model.  Instead we can take any manifold
with an arbitrary number of asymptotically flat ends and cut it off and rescale
it down appropriately to glue it into a collection of small balls.
\end{example}

We now prove some technical lemmas regarding almost unseen sets.

We begin by noting that most directions in terms of the volume
of $S^{n-1}$ behave in an almost isotropic manner when
a space is almost isotropic off of an almost unseen set (recall Definition~\ref{unseen}).

\begin{lemma}\label{est2}
Given $n \in \NN, R>0$, 
and a map $f:(0,R)\times(0,\infty)\to \NN$.  
For all $h >0$ , and for all 
\be
\epsilon \in (0, \min \{\frac{h}{2}, \frac{\pi}{4f(h/2, R+h/2)} \},
\ee
there exists 
\be
\theta(n,f,R,h):= \frac {\pi}{2f(h/2, R+h/2)} %actually very rough estimate here
\ee
such that if $M^n$ is $\epsilon,R$ almost isotropic off 
an $(\epsilon,R)$-almost unseen set
and satisfies the 
$f$ ball packing property then for all $t\in [0,R)$,
\be \label{est3}
F_{q}(\theta, t,t)< h \qquad \forall t\in [0,R), \, \forall \theta<\theta(n,f,R,h).
\ee
\end{lemma}

Note that one cannot expect to control $F$ better than $\epsilon$
because its behavior is defined by $\epsilon$.  Note also
that the bound on $\theta$ does not depend on $\epsilon$.

\Pf 
Suppose on the contrary that
\be 
F_q(\theta,t,t)\ge h.
\ee
Since $F$ is nondecreasing in
$\theta$ and $M$ is almost isotropic this means that 
\be \label{esteq1}
d_M(exp_q(tv), exp_q(tw))>h-\epsilon>h/2,
\ee
whenever $v,w\in S_q$ such that $d_S(v,w)\ge \theta$.  
Recall that
\be
S_q\subset \bigcup_{j=1}^N B_{w_j}(\epsilon_j) \subset S^{n-1} \subset TM_q
\ee
in Definition~\ref{unseen}.  

Now in $S^{n-1}$ with the standard metric $d_S$, there are
at least
\be
N_\theta=\frac{\pi}{(\theta+2\epsilon)}
\ee
disjoint balls of radius $\theta/2 + \epsilon$.
If any of these balls is centered in $S_q$, then it is centered
in a ball $B_{w_j}(\epsilon_j)$, and so it contains a ball
of radius $\theta/2 <\theta/2+\epsilon-\epsilon_j$ centered
in 
\be
B_{w_j}(3\epsilon_j)\setminus B_{w_j}(\epsilon_j)\subset S^{n-1}\setminus S_q.
\ee

Let
$v_1, v_2,... v_{N_{\theta}}$ be the centers of these balls.
Then by (\ref{esteq1}), $B_{exp_q(tv_i)}(h/2)$ are
disjoint as well and contained in $B_q(t+h/2)\subset B_q(R+h/2)$
.  But by the $f$ ball packing property
there are at most $f(h/2,R+h/2)$ disjoint balls of radius
$h/2$ in a ball of radius $R+h/2$.  Thus
$
f(h/2, R+h/2)\ge N_\theta \ge \pi/(\theta+2\epsilon),
$
and
\be
\theta\ge \frac{\pi}{f(h/2, R+h/2)} - 2\epsilon \ge \frac{\pi}{2f(h/2, R+h/2)}. 
\ee
\qed

%----------------------------------------------------------------
\sect{Almost Isotropy and Exponential Maps} \label{beforeonto}

We begin with a definition.

\begin{defn}  \label{WY}\label{defWY}
Let 
\begin{equation}
W_Y:=
\{ y\in Y: \textrm{ there exist } x_i\in W_{M_i}
\textrm{ converging to } y \},
\end{equation}
and let
\begin{equation}
W_\infty:=
\{ y\in Y: \textrm{ there does not exist } x_i\in M_i\setminus{W_{M_i}}
\textrm{ converging to } y \}\subset W_Y.
\end{equation}
\end{defn}

Note that
$W_\infty$
 are the points that cannot be examined using the almost isotropy
properties of the $M_i$.
To prove Theorem~\ref{MainThm} we will show $W_Y$ is discrete
using the fact that the $W_{M_i}$ consist of uniformly disjoint balls 
(\ref{Wballs}), however, in this section we will make no assumption
on the $W_{M_i}$ other than the fact that they are ``almost unseen''.
We prove the following theorem.

\begin{theorem}\label{allbutonto}
Let $M_i$ be a  sequence of Riemannian manifolds possibly with boundary
that are $(\epsilon_i,R)$-almost isotropic off an
$(\epsilon_i, R)$-almost unseen set $W_{M_i}$ which includes
the boundary of $M_i$ if it exists.  Suppose further that
$M_i,p_i$ converge to a complete length space $Y,y$
in the pointed Gromov-Hausdorff sense.  

%Include all key properties from this section:

For all $x\in Y\setminus W_\infty$, there is a continuous
map $exp_x: B_0(R)\subset \RR^n \to B_y(R)\subset Y$ which is a homeomorphism
onto its image $exp_x(B_0(R))$.  Furthermore, the curves
$exp_x(tv)$ for $t\in [0,R]$ are length minimizing.
Furthermore there is isotropy in the sense that there
exist functions $F_x: [0,\pi]\times [0,R)\times [0,R) \to [0,2R)$
satisfying (\ref{abc2a})-(\ref{abc2c}) and (\ref{aboiso}).
\end{theorem} 

Note that (\ref{abc2a})-(\ref{abc2c} are just the natural limits of 
(\ref{abc1a})-(\ref{abc1c}) of Definition~\ref{almiso}.

In general the exponential 
map won't be surjective, as can be seen in the case where $Y$ is
a sphere and a plane joined at a point.  In that case $exp_x$
will only map onto the intersection of $B_x(R)$ with the plane containing
$x$.

We prove this theorem through a series of lemmas.

One of the special properties of Gromov-Hausdorff Convergence
is that if we have a sequence of curves, 
$C_i:[0,L_i]\to M_i$,
parametrized by arclength with $L_i \le L$, then a subsequence has a limit 
$C_\infty:[0,L_\infty]\to Y$ which is a 
curve parametrized by arclength (although $L_\infty=\lim_i L_i$ might
be $0$).  This follows from the generalized Arzela Ascoli Theorem
\cite{GrPet}.
Furthermore if the $C_i$ are length minimizing, so is $C_\infty$.

This allows us to make the following definition.

\begin{defn} \label{Cvi}
Let $v_i\in S^{n-1}\setminus S_{q_i}$ such that
the curves $exp_{q_i}(tv_i)$ for $t\in [0,R]$
converge to a limit curve, then $C_{\{v_i\}}:[0,R]\to Y$ is
their limit curve.
\end{defn}

There is no natural relationship between
the $v_i$ from the different tangent
cones $TM_{q_i}$.  For this reason we fix an identification between
all the $TM_{q_i}$.  Each identification is determined only up
to $SO(n-1)$ but we need to make a choice.
Thus all the $S^{n-1}\setminus S_{q_i}$ can be thought of as subsets
of the same $S^{n-1}$.  By Definition~\ref{unseen} it is easy to see that the 
$S^{n-1}\setminus S_{q_i}$ converge
to this $S^{n-1}$.

\begin{lemma}\label{est1}
Suppose $v\in S^{n-1}$ and $v_i, w_i\in S^{n-1}\setminus S_{q_i}$ are both
sequences converging to $v$, such that
the curves $exp_{q_i}(tv_i)$ converge to a limit curve
$C_{\{v_{i}\}}(t)$, then $exp_{q_{i}}(tw_{i})$ 
also converges to the same limit curve without
having to take a subsequence.  In particular
\be \label{est1a}
\limsup_{i\to\infty}F_{q_i}(\theta_i,t,t)=0\textrm{ if } \theta_i \to 0.
\ee
\end{lemma}

\Pf
By the $(\epsilon_i, R)$-almost isotropy
\be
d_{M_i}(exp_{q_i}(tv_i)), exp_{q_i}(tw_i))<
F_{q_i}(d_S(v_i, w_i), t,t) + \epsilon_i \qquad \forall t\in [0,R].
\ee
Now we know a subsequence of $exp_{q_{i_j}}(tw_{i_j})$ must converge
to a limit curve $C_{\{w_{i_j}\}}$.  We need only show that 
$C_{\{w_{i_j}\}}=C_{\{v_i\}}$.

Using the subsequence and taking the liminf as $i_j\to\infty$
on both sides we get:
\be
d_{Y}( C_{\{v_i\}}(t), C_{\{w_{i_j}\}}(t)) \le
\liminf_{i_j\to\infty} F_{q_i}(d_S(v_i,w_i), t,t) 
\qquad \forall t\in [0,R].
\ee
Now $d_S(v_i, w_i)$ converges to $d_S(v,v)=0$, so we are done if 
we can show (\ref{est1a}) which follows from Lemma~\ref{est2}.
\qed

We can now apply Theorem~\ref{arzasc} to our exponential functions.

\begin{lemma}\label{equicont}
If we assume that $M_i\to Y$ are locally $(\epsilon_i,R)$-almost isotropic off
sets $W_{M_i}$ which are $(\epsilon_i,R)$-almost unseen 
and $\epsilon_i$ converges to $0$, then
we can show that the maps, 
\be
exp_{q_i}:[0,R]\times (S^{n-1}\setminus S_{q_i})\to B_{q_i}(R),
\ee
are uniformly almost equicontinuous for all $q_i\in M_i\setminus W_{M_i}$.
Thus for any $x\in Y\setminus W_\infty$, there is a subsequence
of the $i$ with  a continuous limit map 
$exp_x:[0,R]\times S^{n-1}\to B_y(R)$
such that
$exp_x(0v)=x$  for all $v\in S^{n-1}$ and,
\be \label{a-b}
d_Y(exp_x(av), exp_x(bv))\le |b-a|.
\ee
\end{lemma}

It should be noted that at this stage the limit exponential
map is not necessarily an exponential map in the sense that
$exp_x(tv)$ is a minimizing curve parametrized proportional to
arclength.  Nor is it known to be surjective.  It is also
possible that this exponential map depends on the choice of
the sequence of $q_i\in M_i$ converging to $x\in Y$.
Nevertheless we can set up a local isotropy of sorts using these
exponential maps.

\begin{figure}[htbp]
\includegraphics[height=1.8in ]{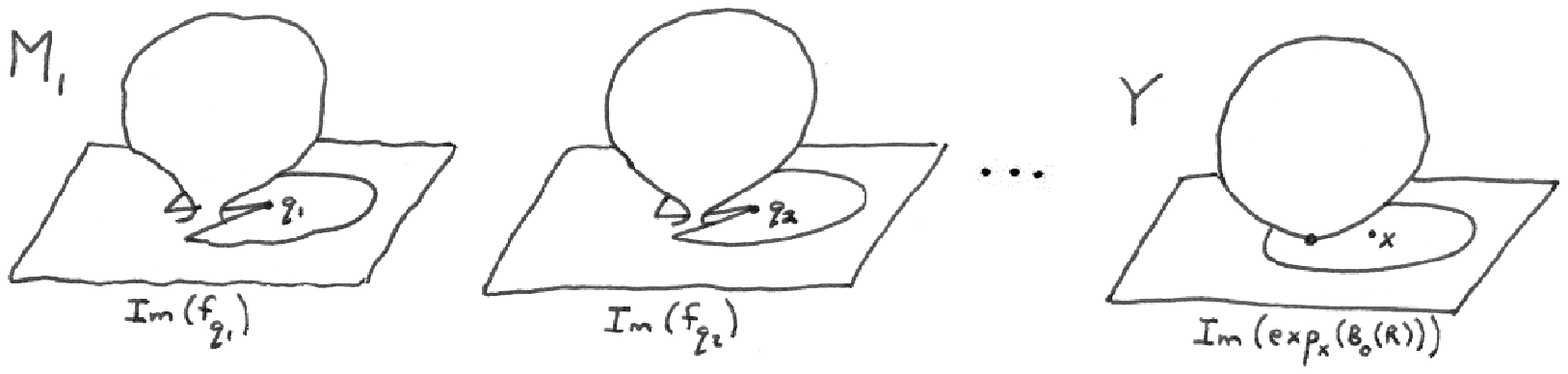}
\caption{Here we see $f_i([0,R)\times (S^{n-1}\setminus S_{p_i})) \subset M_i$,
each looking like a disk with a wedge removed to avoid $W_{M_i}$, 
converging to a limit $exp_x(B_0(R)) \subset Y$,
which is a disk in a plane but not a ball in $Y$.}\label{fignotonto1}
\end{figure}

\Pf
Let $f_i:[0,R)\times (S^{n-1}\setminus S_{p_i}) \to M_i$ be defined
$f_i(s,v)=exp_{p_i}(sv)$.  See Figure~\ref{fignotonto1}.

%\psdraft
%\psfull
%note the diagram is identical to fignotonto2 which could be removed to save room

By the $(\epsilon_i,R)$-almost unseen property, 
$([0,R)\times (S^{n-1}\setminus S_{q_i}), (0,v_0))$
converges to $([0,R]\times S^{n-1}, (0,v_0))$ 
in the pointed Gromov Hausdorff sense.  Furthermore
$f_i(0,v_0)=q_i$ and we are already  given $(M_i,p_i)\to (Y,y)$,
$q_i\to x$ so $(M_i,q_i)$ converges to $(Y,x)$.  So
we need only verify that $f_i$ are uniformly equicontinuous
and can use the functions $F_i$ of the almost isotropy
to do so.

Given any $h>0$, let $\theta<\theta(n,f,R,h/2)$ of 
Lemma~\ref{est2}, so for $i$ sufficiently large that
\be
\epsilon_i <  \min\{\frac{h}{4}, \frac{\pi}{4f(h/4, R+h/4)} \}
\ee
we have
\be 
F_{q}(\theta, t,t)< h/2 \qquad \forall t\in [0,R). 
\ee
So
\begin{eqnarray}
d_{M_i}(f_i(s_1,v_1), f_i(s_2,v_2))&\le &d_{M_i}(f_i(s_1,v_1), f_i(s_1,v_2))
+ d_{M_i}(f_i(s_1,v_2), f_i(s_2,v_2))\\
&<& F_i(d_S(v_1, v_2), s_1, s_1)+(\epsilon_i+|s_1-s_2|) \\
&<& h/2+\epsilon_i+h/2,
\end{eqnarray}
as long as $d_S(v_1, v_2)<\theta$ and $|s_1-s_2|<h/2$.

Thus we are uniformly almost equicontinuous and the rest follows from Lemma~\ref{equicont}.
\qed

\begin{lemma} \label{Fx}
For $Y$ as in Lemma~\ref{equicont} there are 
continuous functions 
$F_x:[0,R]\times S^{n-1} \to R$
defined at each point $x\in Y\setminus W_\infty$, such that
$F$ satisfies conditions (\ref{abc2a})-(\ref{abc2c}) and
\be
F_x(d_S(v,w),t,s)=d_Y(exp_x(tv),exp_x(sw))=
\lim F_{q_i}(d_{S}(v_i,w_i),t,s)
\ee
for any $v,w$ in $S^{n-1}$. 
Here $v_i, w_i \in S^{n-1}\times S_{q_i}$ converge to $v$ and $w$
respectively and $exp_x$ and $F_x$ are defined using the same
sequence of $q_i$.  
\end{lemma}

\Pf
First let $q_i\in M_i\setminus W_{M_i}$ converge to $x$
and let
$v_i\to v$ and $w_i\to w$.  Then by Lemma~\ref{equicont}
$exp_{q_i}(tv_i) \to exp_x(tv)$ and so 
\be
d_Y(exp_x(tv),exp_x(sw))=
\lim_{i\to\infty} d_{M_i}(exp_{q_i}(tv_i),exp_{q_i}(sw_i))=
\lim_{i\to\infty} F_{q_i}(d_{S}(v_i,w_i),t,s).
\ee  
%by squeeze thm
In particular the limit on the right hand side exists.  However this limit
clearly depends only on the angle and the lengths, so
\be
d_Y(exp_x(tv),exp_x(sw))=F_x(d_{S}(v,w),t,s).
\ee
Furthermore, since $F_{q_i}$ satisfy (\ref{abc1a})-(\ref{abc1c}) 
of the definition
of almost isotropy, $F_x$ satisfies (\ref{abc2a})-(\ref{abc2c}).
\qed

\begin{lemma} \label{mingeod}
For $Y$ and $exp_y$ as in Lemma~\ref{equicont},
$d(exp_y(tv), y)=t$ for all $t\in [0,R)$ and so $exp_y(tv)$ is a minimizing
curve parametrized proportional to arclength.
\end{lemma}

Note the geodesics in the $M_i$ were not assumed to be length minimizing but
that the almost isotropy implies that they are almost length minimizing.

\Pf
Fix $v\in S^n$ and $t\le R$ and let $z=exp_y(tv)$.  
Then there exists $v_i\to v$, such that
$z_i=exp_{q_i}(tv_i)\in B_{q_i}(R))$ converge to $z$
by Lemma~\ref{equicont}.  
First note that,
\be
d_Y(y,z)=\lim_{i\to\infty} d_{M_i}(q_i,z_i)\le t.
\ee

On the other hand by the triangle inequality
and (c) in the definition of almost isotropy,
\begin{eqnarray}
d_Y(y,z)&=&\lim_{i\to\infty} d_{M_i}(q_i,z_i)\\
&\ge& \lim_{i\to\infty} d_{M_i}(exp_{q_i}(Rv_i),q_i)-(R-t)\\
&\ge&
\lim_{i\to\infty} F_{q_i}(0,0,R)-(R-t) \,\, =\,\,  R-R+t \,\,  = \,\, t.
\end{eqnarray}
\qed

We also get some interesting properties from the triangle inequality.

\begin{lemma} \label{triF}
$F_x(\theta_1, t_1,s)+F_x(\theta_2, s, t_2) 
\ge F_x(\theta_1+\theta_2, t_1, t_2)$
\end{lemma}

\Pf
Just apply the triangle inequality to $exp_y(t_1v_1)$,
$exp_y(t_2v_2)$ and $exp_y(sw)$ where $d_S(v_i,w)=\theta_i$.
\qed

\begin{corollary}\label{fract}
$F_x(\pi/k,t,t)>t/k$.
\end{corollary}

\Pf
Apply Lemma~\ref{triF} repeatedly, so that
$k F_x(\pi/k,t,t)\ge F_x(\pi,t,t)$
and then apply property (b) of Lemma~\ref{Fx}.
\qed

Although we have defined $exp_x$ as a function of
two variables, a length and a unit vector, 
we know $exp_x(0v)=exp_x(0w)$ for all $v$ and $w$,
so we can consider it as a function of $\RR^n$.

\begin{lemma} \label{1to1}
For fixed $x\in Y$, $exp_x:B_{0}(R)\to B_x(R)$
is a one to one map.
\end{lemma}

\Pf
Suppose not, then there exists $v,w\in S^{n-1}$ and
$t,s \in (0,R]$ such that $exp_x(tv)=exp_x(sw)$.
By Lemma~\ref{mingeod}, $t=s$.  So by Lemma~\ref{Fx}
we have $F_x(\theta, t,t)=0$ for some $\theta>0$.
Since $F_x$ is nondecreasing in its first variable
and there exists $k$ sufficiently large that $\theta>\pi/k$
we have $F_x(\pi/k,t,t)=0$.  This contradicts
Corollary~\ref{fract}.
\qed

Putting these lemmas together we have Theorem~\ref{allbutonto}.

We need only show that $W_\infty$ is discrete to prove that
$Y$ is an exponential length space.  To prove this we need additional 
conditions on the sequence $M_i$.  This can be seen because the $M_i$
could be a pair of planes which are connected by Schwarzschild solutions
at an increasingly dense set of points and still satisfy all the conditions
used in this section.  In the next two sections we show how additional
conditions can be found and satisfied.

%maybe put relatively open here?
%----------------------------------------------------
\sect{Local Surjectivity} \label{sectonto}

In this section we use the condition that the bad sets $W_{M_i}$ which 
are avoided in the definition of the almost isotropy are each contained in 
a union of balls
of decreasing radii that are a uniform distance apart (\ref{Wballs}).

Recall the definitions of $W_\infty \subset W_Y$ in
Definition~\ref{WY} and the exponential map defined in 
Theorem~\ref{allbutonto}.  

\begin{lemma}\label{discrete} 
If $M_i$ converge to a space $Y$ in the Gromov Hausdorff sense and
subsets $W_{M_i}$ satisfy (\ref{Wballs}), then 
$W_Y$ is a countable collection of points $\{y_j\}$ such
that $D_Y(y_j, y_k)\ge 2\bar{R}$ and so $W_\infty$ is an
empty set.
\end{lemma}

\Pf
Given $y_1,y_2 \in W_\infty$.
By the definition of $W_\infty$,
we know there exists $x_i \to y_1$ and $z_i \to y_2$
where $x_i, z_i \in W_{M_i}$.  Since the radius of the
balls in $W_{M_i}$ decreases to $0$, but
$d_{M_i}(x_i,z_i)\to d_Y(y_1,y_2)>0$, eventually
$x_i$ and $z_i$ will be in distinct balls,
and thus $d_{M_i}(x_i,z_i)\ge 2\bar{R}$
and the lemma follows.
\qed

%ADD lemma limvol V(z,s) for s<R IF NEEDED  

\begin{lemma}\label{exponto}
Suppose $M_i\to Y$ satisfy all the conditions of Theorem~\ref{MainThm},
then for all $x\in Y\setminus W_Y$, there is an
$r_x>0$ such that $exp_x:B_0(r_x)\to B_x(r_x)$ is onto.  
\end{lemma}

\begin{figure}[htbp]
\includegraphics[height=1.8in ]{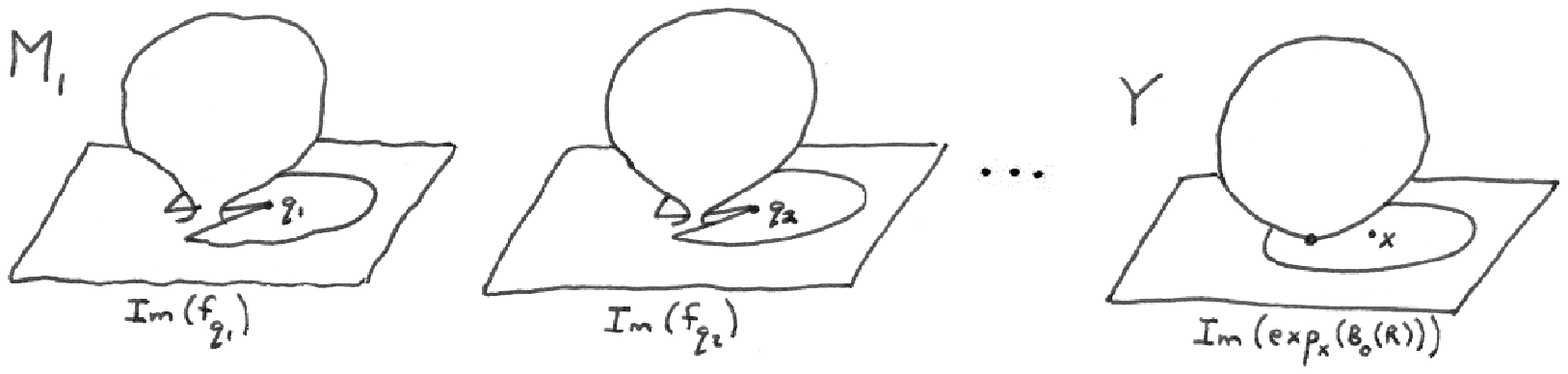}
\caption{Here we see that $exp_x$ does not map $B_0(R)$ onto $B_x(R)$ 
because it's image is a limit of 
$f_i([0,R)\times (S^{n-1}\setminus S_{p_i})) \subset M_i$,
depicted here as disks with wedges removed,
which avoid the $W_{M_i}$
[Lemma~\ref{equicont}].
}\label{fignotonto2}
\end{figure}

Recall that in general it does not map onto
$B_x(R)$ as seen in  Figure~\ref{fignotonto2}.
Note also that $exp_x$ is not shown to map onto any balls about $x$
if $x$ is in $W_Y$.

%\psdraft
%\psfull

\Pf
First by Lemma~\ref{equicont}, we know that for all $x\in Y\setminus W_Y
\subset Y\setminus W_\infty$ we have an exponential function
$exp_x$ which is the limit of some selected subsequence
$exp_{q_i}$ restricted to almost isotropic directions where 
$q_i\in M_i$ converge to $x$.

To prove this lemma we first take $r_x=d_Y(x, W_\infty)/2>0$ since $W_\infty$
is discrete (by Lemma~\ref{discrete}).  
Taking $q_i\in M_i$ as above converging to $x$ 
eventually $d_{M_i}(q_i,W_{M_i})>r_x$ as well,
so $tv,sw \in T_{q_i}$ of Definition~\ref{almiso}
for all $v,w \in S^{n-1}\subset TM_{q_i}$
as long as $t,s\le r_x$.  So we have
\be
|d_{M_i}(exp_{q_i}(tv),exp_{q_i}(sw))-F_{p_i}(d_S(v,w), t,s)|<\epsilon_i.
\ee
Thus we can have almost equicontinuity for
$f_i(t,v)=exp_{q_i}(tv)$, $f_i: [0,r_x)\times S^{n-1}$
and a subsequence converges.  Since we already had $f_i(tv)$
restricted to $[0,r_x)\times S_{q_i}$ converges to
a continuous function $exp_p(tv)$, these limits must agree.
So in fact for any $v_i$ converging to $v\in S^{n-1}$
and $t_i\to t$ in $[0,r_x)$, we have,
\be
exp_x(tv)=\lim_{i\to\infty} exp_{q_i}(t_iv_i).
\ee

Now for any $z\in B_x(r_x)$ there exist 
\be
z_i\in B_{q_i}(r_x)\subset M_i\setminus W_{M_i}
\ee
converging to $z$ and there exist 
\be
v_i\in B_0(r_x)\cap T_{q_i}
\textrm{ such that }
 exp_{q_i}(v_i)=z_i,
\ee 
so a subsequence of the $v_i$ converge
to some $v\in B_0(r_x)$ such that $exp_x(v)=z.$
\qed

Note that the uniformly disjoint ball condition, (\ref{Wballs}),
is required for this lemma to hold as can be seen here.  

\begin{example}
Let $M_i$ be a pair of planes each 
with $i$ balls removed, $B_{(0, 1/k)}(1/(ik^2))$
where $k=1..i$.  Replace the pairs of corresponding balls with
smoothly attached Schwarzschild solutions as in Lemma~\ref{schcut}. 
Then $M_i$ are locally
$(\epsilon_i, 1)$-almost isotropic off $W_{M_i}$ equal to
the collection of edited balls.  These $W_{M_i}$ are 
$(\epsilon_i, 1)$-almost unseen but do not satisfy
(\ref{Wballs}).  The limit space $Y$ is a pair of planes
joined at the points $(0,1/k)$ where $k=1,2,...$.  This set
of points is $W_\infty$.  The
limit exponential map based at $(0,0)$ on one plane does
not map onto any points in the other plane except for those
in $W_\infty$, and so $exp_{(0,0)}$ does not map onto any
balls no matter how small.
\end{example}
  
Note that Lemma~\ref{exponto} and Theorem~\ref{allbutonto}
imply that the limit space $Y$ of $M_i$ satisfying the
combined hypotheses %REF HERE?
is a locally isotropic exponential length space off of $W_Y$
where the exponential maps are chosen depending on 
the $p_i\in M_i$ converging to $x\in Y$ and 
the subsequences.  Later we will
show that in fact no subsequences or choices were necessary
because the limit space will be a space form regardless of the choices
made and the limit space is unique.
 
I conjecture that one could replace (\ref{Wballs})
with a lower bound on Ricci curvature.  Thus far all arguments 
possibly leading to such a statement are lengthy, and so this
question has not been pursued in this paper.

%--------------------------------------------------------------
\sect{Exponential Length Spaces} \label{sectexpl}

This section
focuses on exponential length spaces, $Y$,
off a set of discrete points $W_Y$.
Recall definitions~\ref{explsp} and~\ref{explsp2} from
the Introduction.  
%Recall also that a space form is a complete 
%Riemannian manifold with constant sectional curvature.
To simplify notation, $Y'$ will be a connected component of
$Y\setminus W_Y$ and its closure will be denoted $Cl(Y')$. 
Ultimately we will prove Theorem~\ref{isoexplthm} that $Cl(Y')$
is a space form.

\begin{example} \label{sphpl}
A sphere $\SSS^2$ and a plane $\EE^2$ joined at a point $p$
is an exponential length space off that point.  The exponential
function can be taken to agree with the exponential function
of each space form and is not defined at the point $p$.
Note that the image of $exp_q$ for $q\in \SSS^2$ is contained
in $\SSS^2$, so that in this sense the two components of
the space with $p$ removed don't ``see'' each other.
\end{example}

We now prove some lemmas about exponential length spaces
that are not necessarily locally isotropic.

\begin{lemma}\label{moremin}
If $Y$ is an exponential length space off of a set $W_Y$
then for all $y\in Y\setminus W_Y$,
\be
d_Y(exp_y(sv), exp_y(tw))\ge |s-t| \qquad \forall s,t\in [0,R), v,w\in S^{n-1}.
\ee
\end{lemma}

\Pf
If not, there exists $s>t\in [0,R)$ and $ v,w\in S^{n-1}$ such that
\be
d_Y(exp_y(sv), exp_y(tw))< s-t.
\ee
But then by the triangle inequality,
\be
d_Y(exp_y(sv), y)< s-t +d_Y(y, exp_y(tw))=s-t+t=s
\ee
which contradicts the length minimizing property of the exponential map.
\qed

\begin{lemma}\label{relopen}  
If $Y$ is an exponential length space off a set $W_Y$
then for all $x\in Y\setminus W_Y$,
$exp_x$ is a homeomorphism from $B_0(R_x)$ to its image
$exp_x(B_0(R_x))$.  In particular, if
$q_i,q \in exp_x(B_0(R_x))$ and $q_i$ converge to $q$
then 
\be \label{expinvcont}
\lim_{i\to\infty} exp_x^{-1}(q_i)=exp_x^{-1}(q).
\ee
\end{lemma}

\Pf
Since $exp_x:B_0(R_x) \to exp_x(B_0(R_x))$ is continuous and 1:1
and onto, we need only show that the inverse map is continuous
in the sense described in (\ref{expinvcont}).
This convergence is the
same whether it is defined in the relative topology of
$exp_x(B_0(R_x))$ or on $Y$ itself.  Note that 
there exists an $\epsilon>0$ such that $d_Y(q,x)<R_x-\epsilon$
so eventually the $d(q_i, x)<R_x-\epsilon/2$.
Then $v_i =exp_x^{-1}(q_i)\in B_0(R_x-\epsilon/2)$
have a converging subsequence to some $v_\infty \in B_0(R_x)$.
By the continuity of $exp_x$, 
\be
exp_x(v_\infty)=\lim_{i\to\infty} exp_x(v_i)=\lim_{i\to\infty}q_i=q.
\ee
So $v_\infty=exp_x^{-1}(q)$ and this is true for any subsequence of
the $v_i$.  Thus the limit of the $v_i$ exists and we get (\ref{expinvcont}).
\qed 

Note that the images under $exp_x$ of open sets are {\em relatively open}
in $exp_x(B_0(R_x))$ but are not necessarily open.  In Example~\ref{sphpl} 
if $x\in \SSS^2$ is taken to be within distance $\pi/2$ of the plane, then
$exp_x(R_x)$ is contained completely in $\SSS^2$ and yet it contains
the  point in $W_Y$ which is not in its interior.

By the definition of an exponential
length space, we know that $exp_x$ maps onto $B_x(r_x)$ thus we have 
the following.

\begin{corollary} \label{rx}
If $Y$ is an exponential length space off a set $W_Y$
then for all $x\in Y\setminus W_Y$,
$exp_x$ is a homeomorphism from $B_0(r_x)$ to $B_x(r_x)$.
\end{corollary}

The following Lemma now describes the whole image of the exponential
map proving that the images under $exp_x$ of open sets are indeed
{\em relatively open} in $Y'$.

\begin{lemma} \label{ontoconn}
If $Y$ is an exponential length space off a discrete set $W_Y$,
let $Y'$ be a connected comp of $Y \setminus W_Y$.  Then 
for all $y\in Y'$, $exp_y:B_0(R_y) \to B_y(R_y)\cap (Cl(Y'))$
is onto and is therefore a homeomorphism.
\end{lemma}

\Pf
Let $z \in B_y(R)\cap Y'$.  
By the discreteness of $W_Y$ and connectedness of $Y'$,  
there is a continuous curve $c:[0,1] \to B_y(R)\subset Y'$
which is not necessarily minimizing
from $c(0)=y$ to $c(1)=z$.   Let
\be
T:= c^{-1}(exp_y(B_0(R_y))).
\ee
We need only show that $1 \in T$.

Let $t_0=\sup T $.  So there exists $t_i \in T$ such that $c(t_i)=exp_y(v_i)$
where $|v_i|=d_Y(c(t_i), y)$ converges to $d_Y(c(t_0), y)<R_y$ by the
continuity and location of $c$.  Thus a subsequence of the $v_i$ converges
to $v \in B_0(R_y)$ and, since $exp_y$ is continuous, $exp_y(v)=c(t_0)$ and $t_0\in T$.

\begin{figure}[htbp]
\includegraphics[height=2.5in ]{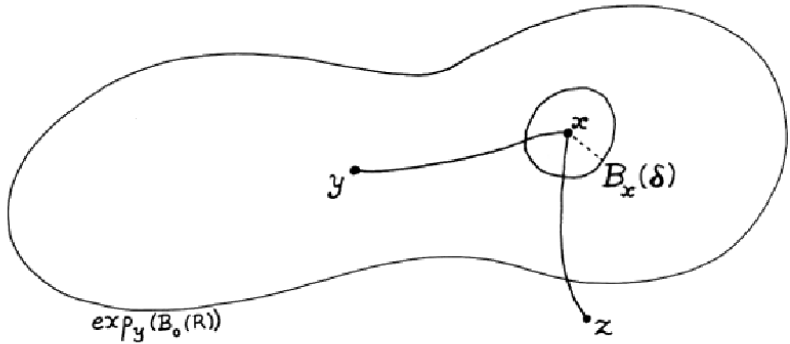}
\caption{Applying the Invariance of Domain Theorem.} \label{fig3}
\end{figure}

Let $x=c(t_0)$ and let $\delta>0$ such that $\delta<min\{r_x, R-d(x,y)\}$,
so $exp_x^{-1}$ is defined on $B_x(\delta)$.  Please consult Figure~\ref{fig3}
while reading this proof.

%\psdraft
%\psfull

Since $exp_y$ is continuous, $exp_y^{-1}(B_x(\delta))$ is an open set in $B_0(R)$ containing $v$,
so there is $\delta'>0$ such that $B_v(\delta')\subset exp_y^{-1}(B_x(\delta)$.

Let $U=exp_x^{-1}(exp_y(B_v(\delta')))\subset \RR^n$.  We claim that $U$ is homeomorphic to $B_v(\delta')\subset \RR^n$.
This would then imply that $U$ is an open set by the Invariance of Domain Theorem
(c.f. \cite{EilSt}) which states that if
$U_i$ are homeomomorphic subsets of $\RR^n$ then
$U_1$ open implies that $U_2$ is open as well.    
Then $exp_x(U)$ is open as well by Corollary~\ref{rx}.  Since 
$c$ is continuous, $c^{-1}(exp_x(U))$ is open, but 
\be
t_0\in c^{-1}(exp_x(U))\subset c^{-1}(exp_y(B_0(R_y)))=T,
\ee
so $t_0=sup T=1$.

Thus it suffices to prove our claim.  To do so, we need only show that 
\be
exp_x^{-1} \circ exp_y: B_v(\delta') \to U
\ee
is continuous and so is its inverse (since we already know it is 1:1 and onto).  
Now, $exp_y$ is continuous by definition and
$exp_x^{-1}$ is continuous by Corollary~\ref{rx} and the fact that $exp_y(B_v(\delta'))\subset B_x(r_x)$.
On the other hand, if $u_i \in U$ converge to $u\in U$ then $exp_x(u_i)$ converges to $exp_x(u)$ by the
continuity of $exp_x$ and since $exp_y^{-1}$ is continuous on $exp_y(B_0(R_y))$ as in (\ref{expinvcont})
we get $exp_y^{-1}(exp_x(u_i))$ converges to $exp_y^{-1}(exp_x(u))$ and we are done.
\qed

The author
would like to thank Prof. Vasquez of CUNY for drawing her
attention to the Invariance of Domain Theorem.

\begin{lemma} \label{univcov}
If $Y$ is an exponential length space off a discrete set $W_y$,
then $Cl(Y')$ has a simply connected universal cover.
Furthermore if $Y$ is locally isotropic then so is the
universal cover of $Cl(Y')$.
\end{lemma}

\Pf
For all $x \in Cl(Y')$ there exists $y\in Y'$ such that $x\in B_y(R)$.
Since all $B_y(R) \cap Cl(Y')$ are simply connected, this means
that $Cl(Y')$ is locally simply connected.  We can then lift the
exponential length structure on these balls isometrically up to
the universal cover.
\qed

Note that the connectedness of $Y'$ is a necessary condition
in the above Lemma as
it is possible that $Y$ could be a bouquet of length spaces 
attached at a point in $W_Y$.  Recall 
from Definition~\ref{explsp} that all $x\in Y\setminus W_Y$ have an
$r_x>0$ such that $exp_x:B_0(r_x)\to B_x(r_x)$ is a homeomorphism
without having to restrict to $Cl(Y')$.

Finally we close with a useful little lemma to deal with the fact
that noncompact exponential length spaces may not be bounded below.

\begin{lemma} \label{R'y}
If $Y$ is an exponential length space off $W_Y$ and then
there exists a positive function $R'_y<R_y/4$ such that
if $x_i\in B_y(R'_y)$ then $R_{x_i}>R'_y$ and $d(x_i,x_j)<R'_{x_i}$.
\end{lemma}

\Pf
Just take $R'_y$ small enough that for all 
$x_i\in B_y(R'_y)$, $R_{x_i} >R_y/2$, so
$d(x_1, x_2)< R_y/4+R_y/4 < R_{x_i}$.
\qed

%---------------------------------------------------
\sect{Locally Isotropic Exponential Length Spaces}\label{sectliso}

Recall the definition of a locally isotropic
exponential length space in Definition~\ref{explsp2}.

We begin with two definitions: the first is classical and the second
is useful for our purposes.

\begin{defn} \label{deflocmin} {\em 
A length space is {\em locally minimizing} if locally
there exist unique minimizing curves between pairs of points.}
\end{defn}

\begin{defn} \label{defulocalmin} {\em 
A length space is {\em uniformly locally minimizing} if 
there exists $R>0$ such that if $d(x,y)<R$ then
there exist unique minimizing curves between $x$ and $y$.}
\end{defn}

Recall that the existence of a minimizing geodesic is a global property
of all complete length spaces by definition.  Here we will show uniqueness.
Note that it is necessary to assume that $Y$ is locally isotropic
to state this lemma.  The standard cone over a circle with an opening
angle $<\pi$ is an exponential length space and it is not locally minimizing.

\begin{lemma}\label{localmin}
Let $Y$ denote an everywhere locally isotropic exponential length space.  
Then for
any length minimizing  curve $c(s)$ from $y$ to $x\in B_y(R)$, 
there exists a unique $v\in S^{n-1}$ such that
$c(s)=exp_y(sv)$.  
\end{lemma}

This has an immediate corollary which follows from Lemma~\ref{R'y}.

\begin{corollary} \label{localmin2}
Let $Y$ denote an everywhere locally isotropic exponential length space.  
Then it is locally minimizing and if $Y$
has a positive exponential radius $R$ then
it is uniformly locally minimizing.
\end{corollary}

Before we prove this lemma we need a technical lemma.

\begin{lemma} \label{Fyinc}
If $Y$ is a locally isotropic exponential length space off a 
set $W_Y$, then for all $y\in Y\setminus W_Y$, $a,b \in (0,R)$,
$\theta \in (0,\pi)$, $F_y(\theta, a,b)>F_y(0,a,b)$.
\end{lemma}

\Pf
Recall that by by (\ref{abc2b}) $F_y(\pi,a,a)>0$.
Given any $\theta>0$, there exists a natural number $k$ such that
$\theta>\pi/(2k)$ so by (\ref{abc2a}), $F_y(\theta, a, b)>F_y(\pi/(2k),a,b)$.
However by the triangle inequality applied to a polygon
of $2k$ points alternatively distances $a$ and $b$ from $y$, we know
\be
k(F_y(\pi/(2k),a,b)+F_y(\pi/(2k),b,a)) \ge F_y(\pi,a,a).
\ee
Thus by the symmetry of $F_y$, we have $F_y(\pi/(2k),a,b)=F_y(\pi/(2k),b,a)>0$.
\qed

We can now prove our Lemma~\ref{localmin}.

\noindent {\bf{Proof  of Lemma~\ref{localmin}:}} 
Let $C:[0,L]\to Y$ be a length minimizing curve from $y$ to $x$.
Then $L=d_Y(x,y)<R$ so $C([0,L])\subset B_y(R)$.
Since $W_Y=\emptyset$ in our hypotheses, we have $r_y=R_y=R$ and Corollary~\ref{rx}
implies that $exp_y^{-1}:B_y(R) \to B_0(R)$ is continuous.  Thus we can define
a continuous map $v(s)=exp_y^{-1}(C(s))$ and we need only show that
$v(s)/|v(s)|$ is constant for  $s>0$.

If not then there exists $s\in (0,L)$ such that $v(s)/|v(s)|\neq v(1)/|v(1)|$
and they have some angle $\theta$ between them.  Thus
using $C$ is parametrized by arclength and Lemma~\ref{Fyinc} we have,

\begin{eqnarray}
L=d_Y(y,x)&=& d_Y(x, C(s))+ d_Y(C(s),y)=d_Y(exp_y(Lv(1)), exp_y(v(s)) +s= \\
&=& F_y(\theta, L, s)+F_y(\theta, s,0) >  F_y(0,L,s) + s = (L-s)+s=L 
\end{eqnarray}  
which is a contradiction.

Thus we can set $v=v(1)/|v(1)|$ and
since we know $exp_y$ is $1:1$ on $B_0(R)$, $v$ must be unique.  

Note  $Y$ is locally minimizing because on any $B_x(R/2)$,
any pair of points is at most $R$ apart.
\qed

We can now prove a significantly stronger technical lemma.

\begin{lemma} \label{solveF}
If $Y$ is a locally isotropic exponential length space off a 
set $W_Y$, then for all $y\in Y\setminus W_Y$, $a,b \in (0,R)$,
we can solve $F_y(\theta, a,b)=d$ uniquely for $\theta$ if a solution
exists unless $a=b=0$.
\end{lemma}

\Pf
Suppose $F_y(\theta_1,a,b)=F_y(\theta_2,a,b)=F_0$ with $\theta_2>\theta_1$.
Without loss of generality because $F_y$ is symmetric in its last two variables
we may assume $b\ge a$.

Then, by (\ref{abc2a}), 
$F_y(\theta,a,b)=F_y(\theta_2,a,b)=F_0$ for all $\theta \in [\theta_1,\theta_2]$.
So
\be \label{solvef1}
d_Y(exp_y(aw), exp_y(bv) )=F_0 \forall v,w \textrm{ s.t. }
d_S(v,w)\in [\theta_1, \theta_2].
\ee

Fix $v_0, w_0\in S^{n-1}$ such that $d_S(v_0,w_0)=\theta_2$
and look at the triangle between the points $y$, $exp_y(aw_0)$
and $exp_y(bv_0)$.  Join the latter two points by a length minimizing
curve $c(t)$ parametrized by arclength
such that $c(0)=exp_y(bv_0)$.  Since $b<R$ we know
$c(t)\in B_y(R)$ for $t<R-b$.  Thus since $exp_y^{-1}:B_y(R)\to B_0(R)$
is a continuous map, there exists continuous curves
$v(t)\in S^{n-1}$ and $b(t)<R$
such that $c(t)=exp_y(b(t)v(t))$ for $t<R-b$.
Now $v(0)=v_0$, so for small $t$ $d_S(v(t),v_0)<\theta_2-\theta_1$
so $d_S(v(t), w_0)\in (\theta_1, \theta_2)$ which combined with
(\ref{solvef1}) implies
that $d_Y(exp_y(bv(t)), exp_y(aw_0))=F_0$.

By the triangle inequality and the fact that $exp_y(b(t)v(t))=c(t)$
running minimally towards $exp_y(aw_0)$, we know
\begin{eqnarray}
F_0&=&d_Y(exp_y(bv(t)), exp_y(aw_0))\\
&\le& d_Y(exp_y(bv(t)), exp_y(b(t)v(t))+ d_Y(exp_y(b(t)v(t)), exp_y(aw_0))\\
&=& |b-b(t)| + (F_0-t).
\end{eqnarray}
Thus $|b(t)-b |\ge t $ for small $t>0$.

However, by Lemma~\ref{localmin}, $exp_y(sv_0)$ is the unique
minimizing curve joining $y$ to $exp_y(bv_0)$ so
\be
b=d_Y(exp_y(bv_0), y)
< d_Y(exp_y(bv_0), exp_y(b(t)v(t))+ d_Y(exp_y(b(t)v(t)), y)= t + b(t).
\ee
Furthermore $exp_y(sv(t))$ is the unique
minimizing curve joining $y$ to $exp_y(b(t)v(t))$ so
\be
b(t)=d_Y(exp_y(b(t)v(t)), y)
< d_Y(exp_y(bv_0), exp_y(b(t)v(t))+ d_Y(exp_y(bv_0), y)= t + b.
\ee
Thus $|b(t)-b|<t$ and we have a contradiction.
 
\qed

\begin{lemma} \label{soisom} 
Suppose $Y$ is a locally isotropic exponential length space off a
discrete set $W_Y$.  
If $y\in Y'$ and $g\in S0(n)$ then there is an isometry
$f_g:B_{y}(R_y)\cap Cl(Y') \to B_{y}(R_y)\cap Cl(Y')$ such that
$f_g(x)=exp_{y}(g(exp_{y}^{-1}(x)))$.
\end{lemma}

\Pf 
By Lemma~\ref{ontoconn}, we know
$exp_{y}:B_0(R_y)\to B_y(R_y)\cap Cl(Y')$ is a homeomorphism
so the inverse is well defined.  We need only verify that
$f_g$ is an isometry.

For any $x_1, x_2\in B_{y}(R_y)$, we know
there exists $s_i\in [0,R_y)$ and $v_i\in S^{n-1}$ such that $exp_y(s_iv_i)=x_i$,
and since $g$ is an isometry on $S^{n-1}$, we have
\begin{eqnarray}
d_Y(f_g(x_1), f_g(x_2))&=& F_{y}(d_S(g(v_1), g(v_2)), s_1, s_2)  \\
& = & F_{y}(d_S(v_1, v_2), s_1, s_2)  \\
& = & d_Y(x_1, x_2).
\end{eqnarray}
Thus we are done.
\qed

\begin{lemma} \label{isotropy} 
Suppose $Y$ is a locally isotropic exponential length space off a
discrete set $W_Y$.  
Given any $y\in Y$, there exists $R=R'_y>0$ such that
if $y_1, y_2\in B_y(R)\cap Y'$ have a length minimizing curve
$\gamma$ running from $y_1$ to $y_2$ of length $d$
such that $\gamma(d/2)\notin W_Y$, then there is an isometry 
$f:B_{\gamma(d/2)}(R)\cap Cl(Y')\to B_{\gamma(d/2)}(R))\cap Cl(Y')$
which fixes $\gamma(d/2)$ and maps $y_1$ to $y_2$ and
$y_2$ to $y_1$.
\end{lemma}

\Pf 
First we set $R'_y>0$ as defined in Lemma~\ref{R'y}.  So
$R_{y_1}, R_{y_2}>R$.

By Lemma~\ref{ontoconn}, we know
$exp_{\gamma(d/2)}:B_0(R)\to B_y(R)\cap Cl(Y')$ is a homeomorphism.
Furthermore $exp_{\gamma(d/2)}^{-1}(\gamma(0))$ and 
$exp_{\gamma(d/2)}^{-1}(\gamma(d))$
are both vectors of length $d/2<R$.  Thus there is an isometry of $\RR^n$ 
fixing $0$ that interchanges these two vectors.  Call it $g$.

So we can define the isometry 
\be
f:B_{\gamma(d/2)}(R)\cap Cl(Y')\to 
B_{\gamma(d/2)}(R))\cap Cl(Y')
\ee
as in Lemma~\ref{soisom},
\be
f(x)=f_g(x)=exp_{\gamma(d/2)}( g( exp_{\gamma(d/2)}^{-1}(x)) )
\ee
By the definition of $g$, $f_g(y_1)=y_2$ and visa versa.
\qed

\begin{coro} \label{Fconst2}%changed!
Suppose $Y$ is a locally isotropic exponential length space off 
a discrete set $W_Y$ and $x,y \in Cl(Y')$, then there is an isometry
\be
f: B_x(R/2)\cap Cl(Y') \to B_y(R/2)\cap Cl(Y')
\ee
where $R=\min_{z\in B_y(2d(x,y)} R'_{z}$.
\end{coro}

Note that $x,y$ must be in the closure of the
same connected component $Y'$ or this is not true as can
be seen when $Y$ is a sphere joined to a plane at a point.

\Pf
First we assume $x,y \in Y'$.
Let $C:[0,L]\to Y'\cap B_y(2d(x,y))$ be a piecewise length 
minimizing curve running from $x$ to $y$.
We only need to show that for all $t$ there is an isometry
$f_t: B_x(R)\cap Cl(Y') \to B_{C(t)}(R)\cap Cl(Y')$ for all $t\in [0,L]$.

We know $f_0$ exists trivially.
Now if $f_s$ exists then for  $s$ near $t$, $f_t$
exists as well using the isometry from Lemma~\ref{isotropy}
to get from $B_{C(s)}(R)\cap Cl(Y')$ to $B_{C(t)}(R)\cap Cl(Y')$.
Furthermore by the standard Arzela Ascoli Theorem, 
if $f_{t_i}$ exist and $t_i\to t$ then a subsequence
converges to a limit isometry with the same
domain and range as the required $f_t$.
So $f_t$ is defined on open and closed intervals and we are done.

Now suppose $x,y \in Cl(Y')$.  Then there exists $x_i \in Y'$ and
$y_i\in Y'$ converging to $x$ and $y$ respectively.  By the
above, we have isometries
$f_i: B_{x_i}(R)\cap Cl(Y') \to B_{y_i}(R)\cap Cl(Y')$.  For
all $r<R$ we can restrict these isometries to the closed
ball $\bar{B}_x(r)$ and we can apply Arzela Ascoli to say that
a subsequence converges to a map 
\be
f_r:\bar{B}_x(r)\cap Cl(Y')\to B_y(R)\cap Cl(Y')
\ee
which preserves distances and is 1:1.
In particular $f_{R/2}$ is an isometry from
$\bar{B}_x(R/2)\cap Y'$ to $\bar{B}_y(R/2)\cap Y'$.
\qed

\begin{lemma} \label{badsexp}
If $Y$ is a locally isotropic exponential length space off
a discrete set $W_Y$
then we can we can define a locally isotropic exponential
length structure everywhere on $Y'$ which is isometric to the original
metric restricted to $Y'$.  This new exponential structure has
$exp_x:B_0(R_x)\to B_x(R_x)\cap Cl(Y')$.
\end{lemma}

We do not claim this extension is unique and clearly the extension
will depend on which connected component of $Y\setminus W_Y$ we are completing.
Note any radius $<R$ will do just as in Lemma~\ref{Fconst2}.

Later in Lemma~\ref{match} we will show to what extent 
isometries preserve length structures.  See also Example~\ref{nomatch} below.
It should also be noted that we have not claimed that $inf_{y\in Y} R_y >0$.
In fact $Y$ could be a Riemannian manifold with constant sectional curvature $-1$
and a cusp end so that $inf_{y\in Y} R_y =0$.

\Pf
For all $x\in Cl(Y')\setminus Y'$, let $R_x=R'_x$ as in Lemma~\ref{R'y} and
let $y\in B_x(R'_x)$.
Then just define 
$exp_x:B_0(R_x/2)\to B_x(R_x/2)$ by taking the isometry 
$f:B_y(R_x/2)\cap Cl(Y')  \to B_x(R_x/2)\cap Cl(Y')$ defined in
Lemma~\ref{Fconst2}.  Then let $exp_x(v)=f(exp_y(v))$
and we are done.
\qed

\begin{lemma} \label{badsexp2}
If $Y$ is a locally isotropic exponential length space 
then we can define a locally isotropic exponential
length structure everywhere on $Y$ which is isometric to the original
metric on $Y$ but has $F_x=F_y$ for all $x,y\in Y$ and $s,t < \min\{R_x, R_y\}$.
\end{lemma}

\Pf
Fix any $y \in Y$ and for any $x\in Y$
define $R$ as in Lemma~\ref{Fconst2} and then define
$exp_x:B_0(R/2)\to B_x(R/2)$ by taking the isometry 
$f:B_y(R/2)  \to B_x(R/2)$ defined in that lemma and let 
$exp_x(v)=f(exp_y(v))$.
\qed

Later in Lemma~\ref{match} we will show this new exponential length structure agrees with
the old exponential length structure in some sense.  However, this is not necessary
at this time.  We will only apply this lemma occasionally and will not in general assume
that $F_y$ is constant.

%---------------------------------------------------------
\sect{Exponential Curves and Locally Minimizing Spaces}\label{sectexpcurve}

In this section we will assume that $Y$ is an exponential length
space everywhere  which is locally minimizing.  So $W_Y=\emptyset$
and $r_x=R_x$ in Definition~\ref{explsp}.
For simplicity, we will take $R_x$ to be small enough both to
satisfy the properties of $R_x$ of the definition of exponential length
space and the $R_x$ of the local minimizing property [Definition~\ref{deflocmin}].  

The following definition should be thought of intuitively 
as the standard differential
equation for a geodesic in a Riemannian manifold adjusted to
make sense in an exponential length space.  We do not yet claim 
that such curves exist and are unique.  

\begin{defn} \label{defexpcurve}{\em
An {\em exponential curve} is a curve $C:[a,b] \to Y$ such 
that for all $t\in [a,b]$, 
there exists $v(t)\in S^{n-1}$ and  satisfying  }
\be \label{diffeq}
C(s)=exp_{C(t)}((s-t)v(t)) 
\quad \forall s\in [t, t+R_t/4]\cap [a,b].
\ee
Here $R_t=R_{c,t}=R'_{c(t)}$ where $R'_x$ is defined in
Lemma~\ref{R'y}. %fancier defn allowing for W_Y existed before december 2002 but unnecessary here
\end{defn}

The following lemma is immediately seen from the definition.

\begin{lemma} \label{overlap}\label{combine}
If $C$ is an exponential curve on $[a,b]$ and on $[b,c]$ and on $[t_1,t_2]$
where $t_1\in (a,b)$ and $t_2\in (b,c)$, then $C$ is exponential on $[a,c]$.
\end{lemma}

\begin{lemma} \label{minexp}
If $Y$ is a %uniformly
locally minimizing exponential length space everywhere 
then all length minimizing curves in $Y$ are exponential curves.
\end{lemma}

Note that if $Y$ has a nonempty $W_Y$ then length minimizing curves 
joining points in distinct
connected components of $Y\setminus W_Y$ are not exponential
as can be seen in the example with the sphere attached to the 
plane at one point [Example~\ref{sphpl}].  

\Pf
Let $c:[a,b]\to Cl(Y')$ be a length minimizing curve.  
For any $t\in [a,b]$ let $t'=\min\{t+R_c, b\}$.
Since $t'-t\le R_b < R'_{c(t)}$ and $Y$ is locally minimizing,
$c([t,t'])$ is a unique length minimizing curve running from $c(t)$ to $c(t')$.
Let $v_t\in S^{n-1}$ be
defined as $exp_{c(t)}^{-1}(c(t'))$.  
Then by the definition of an exponential length space,
we know $exp_{c(t)}((s-t)v(t))$ with $s\in [t,t']\subset [t,t+R_c]$
is also a length minimizing curve from $c(t)$ to $c(t')$.  Thus these curves agree 
and we are done.
\qed

\begin{lemma} \label{extend} 
If $Y$ is a locally minimizing exponential
length space everywhere, the function $exp_x$ can be extended uniquely so that
$exp_x: \RR^n \to Y$, such that
$exp_x(sv)$ is an exponential curve for all $s\in \RR$.
\end{lemma}

%extension can be done even off a set $W_Y$ see p9.7.5.tex

\Pf  Fix $y\in Y$ and $v\in S^{n-1}$.
We treat (\ref{diffeq}) like a differential equation.
We can call the possible solution $C(t)$.

We know $C(t)=exp_y(tv)$ is defined for $t\in [0,R_y]$, so now we must 
extend it.  Clearly if $C(t)$ is defined on an open set it can be
defined on a closed set by extending it continuously.

It suffices to show that if $C(s)$ satisfies
(\ref{diffeq}) for $t\in [0,a]$
then it does for  $t\in [0,a+R_a/8]$ as well.
Although $R_a$ may decrease, we will have proven that $C$
is defined on a right open set, and since it is defined
on a closed set, it is defined  on all of $[0. \infty)$.

Assume $C(s)$ is defined on $[0,a]$. So for all $t\in [0,a]$
\be \label{4.6}
C(s)=exp_{C(t)}((s-t)v_t) \qquad \forall s\in [t,t+R_a/4]\cap [0,a]
\ee
So it is minimizing on $[a-R_a/2,a]$ by the definition of an
exponential length space the fact that $R_a \ge R_{c(a)}$.
%Now choose $s'\in [a-3R/8, a-R/4]$
%such that $C(s')\notin W_Y$. 
Let $s'=a-R_a/4$.  Since 
$ C:[a-R_a/2, a]\to B_{C(s')}(R/2)$ is minimizing
and $R_a \ge R_{c(a-R_a/4)}$ by Lemma~\ref{R'y},
we know that there exists some $w\in S^{n-1}$ such that
\be
C(s)=exp_{C(s')}((s-s')w) \qquad \forall s\in [s',s'+R_{C(s')}] \subset [s',a].
\ee

Extend the definition of $C$ to $[s', a+R_a/4]\subset [s',s'+R_a]$
using this $w$.
\be
C(s) := exp_{C(s')}((s-s')w) \qquad \forall s\in [s',a+R_a/4].
\ee
This extension is a length minimizing curve on $[a-R_a/4,a+R_a/4]$
thus for all $t\in [a, a+R_a/4]$ we know
\be \label{4.9}
C(s)=exp_{C(t)}((s-t)v_t) \qquad \forall s\in [t,a+R_a/4].
\ee
Using (\ref{4.6}) for $t\in [0, a-R_a/4]$ and (\ref{4.9})
for $t\in [a-R_a/4, a+R_a/4]$ we see that $C$ is an exponential
curve on $[0,a+R_a/4]$.  In fact it is slightly better than
exponential since $R_a \ge R_{a+R_a/4}$.  

Since this was true for all $a>0$, $C$ is exponential on $[0,\infty)$.

Note further that $C(t)$ is the only exponential curve which
agrees with $exp_y(tv)$ for $t<R$, since at any point $a$ where they
might split, we are forced to have both satisfy (\ref{4.9})
with the same $v_t$ at $t=a-R/4$. So $exp_y(tv)$ has been extended 
uniquely to all $t\in [0,\infty)$ using this solution $C(t)$.
\qed

We now show this extended exponential map 
is continuous.  It is not 1:1
as can be seen when $Y=\SSS^n$ or $Y=\TT^n$.

\begin{lemma} \label{extcont}
If Y is a locally minimizing exponential length space everywhere
the extended exponential map based at any fixed point is continuous.
\end{lemma}

\Pf
Suppose $t_i \to t$ and $v_i \to v \in S^{n-1}$, we need to show
$exp_y(t_iv_i)$ converges to $exp_y(tv)$. Clearly we
need only show $exp_y(tv_i)$ converges to $exp_y(tv)$
since $|t_i-t|\to 0$, $exp_y(tv_i)$ is parametrized
proportional to arclength and the triangle inequality holds.

By Arzela Ascoli Theorem a subsequence of $c_i(t)=exp_y(tv_i)$ converges to 
a curve $C(t)$ which is parametrized by arclength.  
Well we know $R_t=R_{c_i,t}=\min_{s\in [a,t]} R'_{c_i(s)}/2$ where $R'_x$ is 
defined in Lemma~\ref{R'y} so that $R'_x<R_z$ for all $z$ near $x$.
In particular 
\be
R_{c_i,t} > R_{min,t}=\min_{x\in B_y(t+R_{max})} R_x
\ee
where $R_{max}=\max_{x\in B_y(t)} R'_x$.
So each $c_i$ is a minimizing curve on intervals of length $R_{min,t}$ in $[0,t]$.
Since this is uniform in $i$, the same holds for $C$.

Since $C$ is minimizing on intervals then it must be exponential on those
intervals by Lemma~\ref{minexp}.  Thus it must be an exponential curve by
Lemma~\ref{overlap}.  Since it must agree with $exp_y(tv)$ for small $t$,
it must be its unique extension by Lemma~\ref{extend}.  Thus the
$exp_p(tv_i)$ must have converged without taking a subsequence and we are done.
\qed

To show that the exponential map is open in Riemannian manifolds it is
necessary to avoid conjugate points.  So we need to make a similar
argument in this case.  We can study conjugate points in any
space with an extended exponential map so we will do so in
the following separate section.

%----------------------------------------------------------------------
\sect{Extended Exponential Length Spaces} \label{sectext}

In this section we generalize the properties of the exponential
map and its relationship with conjugate points.  Recall that in Riemannian manifolds
a conjugate point $y$ occurs at $exp_y(t_0v)$ iff $d(exp_y)_(t_0v)$ is
not invertible.  The Implicit Function Theorem is used to show that
a lack of conjugate points on a ball implies that $exp_y$ is a local 
diffeomorphism on that ball.
Here we have no differentiability, but we can use the 
definition of a conjugate point which refers only to length minimizing
curves and we can obtain a local homeomorphism using the Invariance of Domain Theorem.
  
We begin with a definition.

\begin{defn} \label{extexplsp}{\em
A complete length space $Y$ is an {\em extended exponential
length space} if there exists $n\in \NN$
for all $y\in Y$ there exists a map
$exp_y:\RR^n \to Y$ which is continuous 
and there exists a continuous function
$R_y>0$ such that $exp_y:B_0(R_y)\to B_y(R_y)$
is a homeomorphism.  

We also assume that for all $t>s>0$, for all $y\in Y$ and $v\in S^{n-1}$
there exists $w\in S^{n-1}$ satisfying
\be
exp_y(tv)=exp_{exp_y(sv)}((t-s)w)
\ee
and all length minimizing curves are of the form $exp_y(tv)$ with $t\ge 0$.}
\end{defn}

We have already proven that uniformly locally length minimizing exponential
length spaces are extended exponential length spaces in 
Lemmas~\ref{minexp}, \ref{extend} and~\ref{extcont}.  
Conversely, since
the exponential maps in an extended exponential length space
are assumed to be invertible up to some radius $R_y>0$
and since here we assume length minimizing curves are exponential,
we know that extended exponential length spaces are locally minimizing.  

\begin{example} \label{nomatch}
Note that the definition of an extended exponential length space
works only in the positive direction.  We want $exp_y(tv)$
to be an exponential map when $t$ runs from $-1$ to $1$ but this is not
necessarily the case.  For example, one could define an extended
exponential length structure on $\EE^2$ where 
\be
exp_0(t(cos(\theta), sin(\theta))= (tcos (\theta^2/\pi), t sin(\theta^2/\pi))
\ee
for $\theta \in [0, \pi]$ and 
\be
exp_0(t(cos(\theta), sin(\theta))= (tcos (\theta), t sin(\theta))
\ee
for $\theta \in [-\pi, 0]$ and $exp_0(tv)$ would have a corner at $0$.
\end{example}

\begin{lemma}\label{extonto}
If $Y$ is an extended exponential length space then for all $y\in Y$
$exp_y:\RR^n \to Y$ is surjective.
\end{lemma}

\Pf
For all $x\in Y$, there is a length minimizing curve from $y$ to $x$
by the definition of a complete length space.  By Definition~\ref{extexplsp}
that curve must be exponential and have the form $exp_y(tv)$.
\qed

Recall that a {\em cut point} of $y$ has two distinct length
minimizing curves joining it to $y$.  

\begin{lemma} \label{mincut} 
In an  exponential length space $Y$.
If $exp_y(tv)$ is length minimizing on $[0,L]$ then
it has no cut points before $L$.
\end{lemma}

\Pf
Let $x=exp_y(Lv)$.
Suppose that $exp_y(tv)$ has a cut point at $t_0 \in (0,L)$.  
Then there exists two distinct length minimizing curves from 
$x$ to $y$ which both agree with $exp_y((L-s)v)$
for $s\in [0,L-t_0]$ and then diverge.  Since length minimizing 
curves are exponential curves, they cannot diverge,
so they agree everywhere and there is no cut point.
\qed

Now we make a definition of conjugate point for extended exponential
spaces which does not agree exactly with the definition in Riemannian
geometry but is the appropriate extension for our purposes.

 \begin{defn} \label{conjpt} {\em
In an exponential length space, 
an exponential curve $exp_y(tv)$ has a {\em conjugate point} 
$exp_y(t_0v)$ at $t_0>0$ 
if there exists  $v_i\neq w_i$ both converging to $v$ 
and $t_i, s_i \to t_0$
such that $exp_y(t_iv_i)=exp_y(s_iw_i)$.}
\end{defn} 

Recall that in a Riemannian manifold, not all conjugate points take this
form, but any point which has this property is a conjugate point 
(c.f. \cite{doC}).

As in Riemannian manifolds, some conjugate points
are also cut points.  It is easy to see that the first conjugate point
along an exponential curve must have $t_0>R_y$ because $exp_y:B_0(R_y)\to B_y(R_y)$
is one to one. 

We now follow with lemmas extending standard theory of conjugate points
from Riemannian manifolds to this setting. 

\begin{lemma} \label{conjmin}
In an extended exponential length space $Y$:
if $exp_y(tv)$ is an exponential curve with no conjugate
or cut points before $L$, then it is length minimizing on $[0,L]$.
\end{lemma}

\Pf  
Suppose $d_y(exp_y(Lv),y)<L$, then by continuity, for $t$ near $L$
$d_y(exp_y(tv),y)<t$.  Let $t_0=\inf\{t:d(exp_y(tv),y)<t\}\ge R_y$.
Let $s_i$ decrease to $t_0$, then there are 
length minimizing curves $exp_y(tv_i)$ running from $y$ to $exp_y(s_iv)$
of length $t_i<s_i$.
If a subsequence of $v_i$ converges to $v$ then we have a conjugate point.
  
Otherwise a subsequence must converge to some $w$ giving a length minimizing
curve $exp_y(tw)$ from $y$ to $exp_y(t_0v)$.  The latter must then be a cut point.
\qed
 
%TRied CONVERSE see p9.12.06.tex, should check if a counter example, unnecessary here

\begin{lemma} \label{local1to1}
Suppose $Y$ is an extended exponential length space.
If $y\in Y$ has no conjugate points before $t_0>0$ then
$exp_y: B_0(t_0) \to B_y(t_0)$ is locally one-to-one.
\end{lemma}

Note that it is clearly not actually one-to-one as can be seen in 
the cylinder.

\Pf
We need only show that for all $s_0v\in B_0(t_0)$ 
there exists a neighborhood $U\subset B_0(t_0)$ of $s_0v$ 
such that $exp_y:U\to exp_y(U)$ is 1:1.  
If not, $\exists (t_iv_i) \neq (s_iw_i)$ both
converging to $s_0v$ s.t. $exp_y(t_iv_i)=exp_y(s_iw_i)$.
Then $v_i, w_i\to v$ so there is a conjugate point at $s_0<t_0$
unless $v_i=w_i$.  However if $v_i=w_i$ then there are increasingly
small exponential loops converging on the point $exp_y(t_0v)$:
\be
exp_y(t_iv_i)=exp_y(s_iv_i) \textrm{ with } s_i \neq t_i.  
\ee
However if we take $r=R_{exp_y(t_0v)}/2$ then for $i$ sufficiently
large we have $R_{exp_y(t_iv_i)}>r$ and so it cannot
have a loop shorter than $r$.  Taking $i$ possibly larger we get 
$|s_i-t_i|<r$ and a contradiction.
\qed

\begin{lemma} \label{extopen}
Suppose $Y$ is an extended exponential length space.
If $y\in Y$ has no conjugate points before $t_0>0$ then
$exp_y: B_0(t_0) \to B_y(t_0)$  is open.
\end{lemma}

\Pf   
Since there are no conjugate points we can apply Lemma~\ref{local1to1} 
so $\forall v\in B_0(R') \exists \epsilon>0$ s.t. 
$exp_y: B_v(\epsilon) \to im(B_v(\epsilon))$ is 1:1 and continuous.

I claim it is also open.  So I first show
$exp_y^{-1}:im(B_v(\epsilon/2)) \to B_v(\epsilon/2)$ is continuous.
Let $x_i\in im(B_v(\epsilon/2))$ and $x_i\to x\in im(B_v(\epsilon/2))$
Then $\exists v_i \in s.t. exp_y(v_i)=x_i$, and a subsequence of
$v_i$ converging 
to some $v_\infty\in B_v(\epsilon)$.  By continuity of $exp_y$,
$x_i=exp_y(v_i)$ converges to $exp_y(v_\infty)$.  
Thus $exp_y(v_\infty)=x$, but
$exp_y^{-1}$ is unique so all subsequences of the $v_i$ must converge to
the same $v_\infty$, so in fact $v_i$ themselves must converge to $v_\infty$
and $exp_y^{-1}$ is continuous.

\begin{figure}[htbp]
\includegraphics[height=2.5in ]{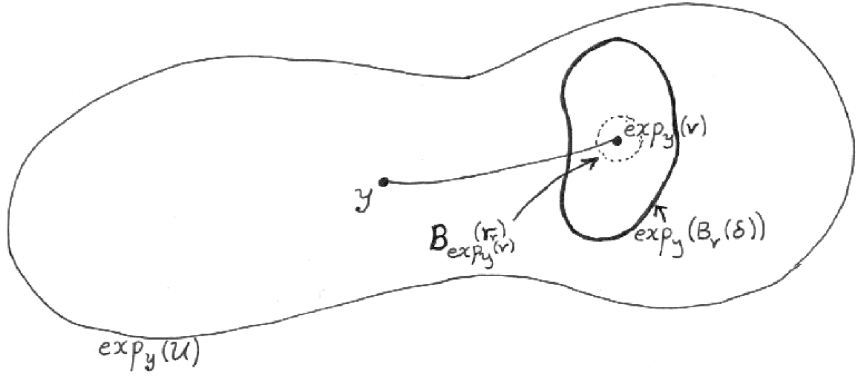}
\caption{ Setting up homeomorphisms from $ B_v(\delta)$ to $ exp_y(B_v(\delta))$
     to $ exp^{-1}_{exp_y(v)}(exp_y(B_v(\delta)))$. } \label{fig4}
\end{figure}

Now we prove that $exp_y$ is open.  Let $U\in B_0(R')$
be any open set,  We must show $exp_y(U)$ is open. 
That is for any $v\in U$, we need to find an $r_v$ such 
that $B_{exp_y(v)}(r_v)\subset exp_y(U)$.  Please consult Figure~\ref{fig4}
while reading this proof.

%\psdraft
%\psfull

Now take any $\delta>0$ sufficiently
 small that $B_v(\delta)\subset U$
and $\delta<\epsilon$.  Thus $exp_y(B_v(\delta))$ is homeomorphic
to $B_v(\delta)$ and is relatively open as a subset of $U$.  
For $\delta$ sufficiently small 
\be \label{IDT1}
exp_y(B_v(\delta))\subset B_{exp_y(v)}(r_{exp_y(v)}).
\ee
Now
\be \label{IDT2}
exp_{exp_y(v)}^{-1}(exp_y(B_v(\delta))) \subset B_0(r_{exp_y(v)}) \subset \RR^n
\ee
is homeomorphic to the open set $B_v(\delta)\subset \RR^n$, so it must be open
by the Invariance of Domain Theorem \cite{EilSt}.

So using $exp_{exp_y(v)}$ as a homeomorphism, we map the sets
in (\ref{IDT2}) back to the corresponding sets in (\ref{IDT1})
and  conclude that $exp_y(B_v(\delta)$ is a relatively open set
in $ B_{exp_y(v)}(r_{exp_y(v)})$, so it is open.
Thus there exists $r_v>0$ such that 
$B_{exp_y(v)}(r_v)\subset exp_y(B_v(\delta)\subset exp_y(U)$
and we are done.  
\qed

\begin{theorem} \label{conjhomeom}
Suppose $Y$ is an
extended exponential length space with no conjugate points about a given
point $y$ before $t=t_0$, then the exponential map $exp_{\tilde{y}}$ is a
local homeomorphism from $B_0(t_0)$ onto $B_y(t_0)$.
So if $Y$ is simply connected it is a homeomorphism.
\end{theorem}

\Pf
The exponential map is continuous, open, locally 1:1 and onto by
Lemmas~\ref{extopen}, \ref{local1to1} and~\ref{extonto}, so
it is a local homeomorphism, thus when applied to a simply
connected space it is a homeomorphism.
\qed

%-------------------------------------------------------------
%check from here on

\sect{Local Isotropy and Conjugate Points} \label{secttube}

%I was fixing tubular neighborhood radius here!!!!

We now return to locally isotropic exponential length spaces
$Y$ with $W_Y$ empty.
Our goal is to show that the universal covers of such spaces
must be $\SSS^n$, $\HH^n$ or $\EE^n$.  Recall that when proving that
complete simply connected Riemannian manifolds with constant sectional
curvature are limited to these three cases, one first shows that
either there are no conjugate points or there is exactly one
conjugate point and thus the exponential map can be inverted
to get either a map from the Riemannian manifold to $\RR^n$ or to $S^n$.

Here we complete this first step [Lemma~\ref{allconj}].
To show that all exponential curves have conjugate points
at the same locations, we first construct isometries
of balls along the exponential curves [Lemma~\ref{ballext}]
and then study the behavior
of exponential curves which are located
near each other [Lemma~\ref{ctogether}].   
(c.f. \cite{doC}).

Recall that we've proven that all locally isotropic exponential length
spaces have extended exponential length structures.
 Recall also that the definition of
an exponential curve is only in a positive direction so that
in general $exp_p(tv)$ is not an exponential curve through $t=0$ [Example~\ref{nomatch}].

\begin{lemma}\label{pi}
For all $y$ in a locally isotropic exponential length space $Y$,
the curve $c(t)=exp_y(tv)$ is exponential for all $t\in \RR$.
\end{lemma}

\Pf
We need only show there exists $\epsilon>0$ such that 
$d_Y(exp_y(-\epsilon v), exp_y(\epsilon v))=2\epsilon$
since all length minimizing curves are exponential
and then we'd be able to apply Lemma~\ref{combine}.

Now
\be
d_Y(exp_y(\epsilon v), exp_y(-\epsilon v))=F_y(\pi, \epsilon, \epsilon) >
F_y(\theta, \epsilon, \epsilon)
\ee
for all $\theta<\pi$ by Lemma~\ref{solveF} combined with (\ref{abc2a}) of
the definition of locally isotropic.    Thus we have a one point set:
\be \label{piunique}
\partial B_y(\epsilon) \cap \partial B_{exp_y(\epsilon v)}(2\epsilon)=\{exp_y(-\epsilon v)\}
\ee

Take $\epsilon<R_y/3$ sufficiently small that if $x \in \bar{B}_y(\epsilon)$
we know $R_x>2R_y/3$.  Set $x=exp_y(\epsilon v)$.
Then there exists $w\in S^{n-1}$ such that $exp_x(tw)=exp_y((\\epsilon-t)v)$
because $exp_y((\epsilon-t)v)$ is minimizing on $(0,\epsilon)\subset [0,R_y]$.
On the other hand $exp_x(tw)$ is minimizing on $[0,2\epsilon]\subset [0,R_x]$
so $d_Y(exp_x(2\epsilon w),x)=2\epsilon$ and by (\ref{piunique})
$exp_x(2\epsilon w)=exp_y(-\epsilon v)$ and we are done.
\qed

\begin{lemma} \label{ballext}
Let $Y$ denote an everywhere locally isotropic exponential length space  
with isotropy radius $R$, where $R$ may be infinity.

Suppose $f:A \to B $ is an isometry
and $A\subset B_x(R/2)$ and $B\subset B_{f(x)}(R/2)$.
Then we can extend the map $f$ so that
$f:B_x(R/2) \to B_{f(x)}(R/2)$ is an isometry.
\end{lemma}

We do not claim this extension is unique.

\Pf
First of all $Y$ can be temporarily given a locally
isotropic exponential length structure such that
$F_x=F_y$ for all $x,y \in Y$ by Lemma~\ref{badsexp2}.

Let $y=f(x)$, $A'=exp_x^{-1}(A) \subset B_0(R)$
and $B'=exp_y^{-1}(B) \subset B_0(R)$.
Let $df_x: A' \to B'$ be defined as
\be \label{df1}
df_x(v)= exp_y^{-1}( f (exp_x(v)) ).
\ee

We claim $df_x$ is in $S0(n)$.  Clearly $df_x(0)=0$
and $|df_x(v)|=|v|$ for all $v\in A'$
by the choice of $x$ and $y$.  
Thus for all $v,w \in A'$ we have
\be
|df_x(v)-df_x(w)|^2=|v|^2+|w|^2 -2|v||w|cos (\theta),  
\ee
where $\theta=d_S(df_x(v), df_x(w))$.
Now using $F_x=F_y$ we have 
\begin{eqnarray}
F_x(\theta, |v|, |w|)&=& d_Y(exp_x(df_x(v)), exp_x(df_x(w))) 
=d_Y( exp_y(v), exp_y(w) ) \\
& =& F_y(d_S(v,w), |v|, |w|) = F_x(d_S(v,w), |v|, |w|).
\end{eqnarray}
This implies $\theta=d_S(v,w)$ by Lemma~\ref{solveF}, so
\be
|df_x(v)-df_x(w)|^2=|v|^2+|w|^2 -2|v||w|cos (\theta)=|v-w|  
\ee
and we have our claim.  

Now since $df_x$ is in $SO(n)$, although it may not
completely be determined depending on the size of $A'$, we can extend $df_x$
to an element of $S0(n)$ mapping $B_0(R)\to B_0(R)$.
Note that $df \cdot  g$ is also a possible extension of $df$
as long as $g$ is in the subgroup of $S0(n)$ preserving $A'$.
 
We can extend the definition of $f$ as
$f(z)=exp_y(df_x(exp_x^{-1}(z)))$ which agrees with $f$ on $A$.  Then
given $x_i\in B_x(R)$ let $s_iv_i=exp_x^{-1}(x_i)$ with $|v_i|=1$
and we have
\begin{eqnarray}
d_Y(f(x_1), f(x_2))&=& F_y(d_S(df_x(v_1), df_x(v_2)), s_1, s_2)  \\
& = & F_y(d_S(v_1, v_2), s_1, s_2) \\
& = & F_x(d_S(v_1, v_2), s_1, s_2)   \\
& = & d_Y(x_1, x_2).
\end{eqnarray}
so we have an isometry.
\qed

\begin{coro} \label{thinsector}
Suppose $Y$ is a locally isotropic exponential length space everywhere.
Given any $L>0$ and $r\in (0,R/4)$
there exists a sector defined by some $\theta_{L,x,r}>0$ such that
\be \label{thinsect1}
d_Y(exp_{x}(tw), exp_x(tv)) <r \qquad \forall t\in [0,L]
\ee
for all $v, w\in S^{n-1}$ such that $d_S(v,w)<\theta_{L,x,v}$.
\end{coro}

\Pf 
If not then 
there exist $w_i, v_i \in S^{n-1}$, $d_S(v_i,w_i)\to 0$ and $t_i \in [0,L]$
such that
\be \label{thinsect2}
d_Y(exp_{x}(t_iw_i), exp_x(t_iv_i)) \ge r.
\ee
Since we know there exist converging subsequences of $t_iw_i$ and of $t_iv_i$,
and they must converge to the same $tv$, and the exponential map is continuous,
we get a contradiction.
\qed

\begin{lemma} \label{ctogether}
Suppose $Y$ is a locally isotropic length space everywhere with isotropy radius $R$
and $N \in \NN$.  Then we can extend the isotropic behavior in thin sectors:
given
\be
d_S(v_1,w_1)=d_S(v_2,w_2)
=\theta <\min\{\theta_{(N+1)R/4,x_1, R/4}, \theta_{(N+1)R/4,x_2, R/4}\}, 
\ee
and
\be
d_S(v_1,\bar{w}_1)=d_S(v_2,\bar{w}_2)
=\bar{\theta} < \min\{\theta_{(N+1)R/4,x_1, R/4}, \theta_{(N+1)R/4,x_2, R/4} \}, 
\ee
and $d_S(w_1, \bar{w}_1)=d_S(w_2, \bar{w}_2)$,
then
\begin{eqnarray} 
d_Y(exp_{x_1}(tv_1), exp_{x_1}(sw_1)) & = & d_Y(exp_{x_2}(tv_2), exp_{x_2}(sw_2)) \label{ctog1}\\ 
d_Y(exp_{x_1}(tv_1), exp_{x_1}(\bar{s}\bar{w}_1)) & = & d_Y(exp_{x_2}(tv_2), exp_{x_2}(\bar{s}\bar{w}_2))\label{ctog1a} \\ 
d_Y(exp_{x_1}(\bar{s}\bar{w}_1), exp_{x_1}(sw_1)) & = & d_Y(exp_{x_2}(\bar{s}\bar{w}_2), exp_{x_2}(sw_2))  \label{ctog1b}
\end{eqnarray}

whenever $t,s, \bar{s} \in [0,NR/4]$ with $|t-s|<R/4$ and $|t-\bar{s}|<R/4$.
\end{lemma}

\Pf
By our choice of $\theta$, for all $k=0..N$ we know that $exp_x([0,kR/4]w_i)\subset {T}_{exp_x([0,kR/4]v}(R/(4k))$
and $exp_x([0,kR/4]\bar{w}_i)\subset {T}_{exp_x([0,kR/4]v}(R/(4k))$.
We will prove (\ref{ctog1})-(\ref{ctog1b}) inductively with $t,s, \bar{s} \in [0.kR/4]$ and k increasing to $N$.

We start with $k=0$.  Choose $g\in S0(n)$ which maps $v_1$ to $v_2$, $w_1$ to $w_2$ and $\bar{w}_1$ to $\bar{w}_2$,  then
by Lemma~\ref{isotropy}, there is an isometry $f_g: B_{x_1}(R) \to B_{x_2}(R)$
which maps $exp_{x_1}([0,R]v_1)$ to $exp_{x_2}([0,R]v_2)$,
$exp_{x_1}([0,R]w_1)$ to $exp_{x_1}([0,R]w_2)$  and $exp_{x_1}([0,R]\bar{w}_1)$ to $exp_{x_1}([0,R]\bar{w}_2)$ 
and we have (\ref{ctog1})-(\ref{ctog1b}) for $t,s \in [0,R/2]$.

Suppose we have (\ref{ctog1})-(\ref{ctog1b}) for $t,s, \bar{s}\in [0,kR/4]$
Set $y_1=exp_{x_1}((kR/4)v_1)$
and $y_2=exp_{x_2}((kR/4)v_2)$ and let 
\be
A=B_{y_1}(R/2)\cap (exp_{x_1}([0, kR/4]v_1)\cup (exp_{x_1}([0, kR/4]w_1)\cup (exp_{x_1}([0, kR/4]\bar{w}_1) ).
\ee
Let
\be
B=B_{y_2}(R/2)\cap (exp_{x_2}([0, kR/4]v_2)\cup (exp_{x_2}([0, kR/4]w_2)\cup (exp_{x_2}([0, kR/4]\bar{w}_2) ).
\ee

By the fact that (\ref{ctog1})-(\ref{ctog1b}) holds for $t,s, \bar{s} \in [0,kR/4]$ 
and all other distances in $A$ and $B$ are determined by the fact that the exponential curves are length minimizing,
we can define
an isometry $f:A \to B$ such that 
$f(exp_{x_1}(tv_1))=exp_{x_2}(tv_2)$,$f(exp_{x_1}(sw_1))=exp_{x_2}(sw_2)$
and 
\be
f(exp_{x_1}(s\bar{w}_1))=exp_{x_2}(s\bar{w}w_2).
\ee
By Lemma~\ref{ballext} we can extend
this isometry to an isometry $f:B_{y_1}(R/2)\to B_{y_2}(R/2)$.  

Since  by the choice of $\theta$,
\be
exp_{x_1}([kR/4-R/4,kR/4+R/4]w_1)\subset B_{y_1}(R/2)
\ee
and,  by the choice of $\bar{\theta}$, 
\be
exp_{x_1}([kR/4-R/4,kR/4+R/4]\bar{w}_1)\subset B_{y_1}(R/2)
\ee
these exponential curves restricted to $[kR/4-R/4, kR/4]$ are in $A$.  Since isometries
map exponential curves to exponential curves we know $f$ must map these
entire segments to the corresponding exponential curves extended as well.  That is
 $f(exp_{x_1}(tv_1))=exp_{x_2}(tv_2)$,
$f(exp_{x_1}(sw_1))=exp_{x_2}(sw_2)$ and $f(exp_{x_1}(s\bar{w}_1))=exp_{x_2}(s\bar{w}_2)$ 
for $t,s\in ((k-1)R/4, (k+1)R/4)$
and we have (\ref{ctog1})-(\ref{ctog1b}) for $t,s\in ((k-1)R/4, (k+1)R/4)$.

Thus by induction we have (\ref{ctog1})-(\ref{ctog1b}) for $t,s, \bar{s}\in ((k-1)R/4, (k+1)R/4)$ where $k=0..N$.
Given any $t,s\in [0,NR/4]$ such that $|t-s|<R/4$ we set 
$k=[4\min_{t,s}/R]$ and  similarly for $\bar{s}$.
\qed

\begin{lemma} \label{minconj}
In a locally isotropic exponential length space $Y$,
if $exp_y(tv)$ is length minimizing on $[0,L]$ then
it has no conjugate points before $L$.
\end{lemma}

\Pf
Suppose that $exp_y(tv)$ has a conjugate point in  $(0,L)$.
So there exists   $v_i\neq w_i$ both converging to $v$ 
and $t_i, s_i \to t_0$
such that $exp_y(t_iv_i)=exp_y(s_iw_i)$.  For $i$
sufficiently large,  
\be
d_S(v_i,w_i)<\min\{\theta_{2L,y, R/4}\} \textrm{ and } |s_i-t_i|<\min\{R/4, L-t_0\}.
\ee
Let $\bar{v}_i$ be chosen such that $d_S(\bar{v}_i, v)=d_S(v_i,w_i)$,
so by Lemma~\ref{ctogether}, we get
\be \label{tisi}
d_Y(exp_{y}(t_iv), exp_{y}(s_i\bar{v}_i)) =d_Y(exp_{y}(t_iv_i), exp_{y}(s_iw_i)) =  0,
\ee
and
\be\label{siti}
d_Y(exp_{y}(s_iv), exp_{y}(t_i\bar{v}_i)) =d_Y(exp_{y}(t_iv_i), exp_{y}(s_iw_i)) =  0.
\ee
Since $exp_y(tv)$ is length minimizing curve up to $L$, and $t_i, s_i<L$,
(\ref{tisi}) implies that 
\be
s_i \le d_Y(y, exp_y(t_iv_i))=t_i
\ee
 while
(\ref{siti}) implies that $t_i \le d_Y(y, exp_y(s_iv_i))=s_i$.
Thus $s_i=t_i$ and both $exp_y(t_i\bar{v}_i)$ and $exp_y(t_iv)$ are 
distinct length minimizing curves which gives us a cut point before $L$.
This is impossible on a length minimizing curve by Lemma~\ref{mincut}.
\qed

\begin{lemma} \label{allconj}
Let $Y$ denote an everywhere locally isotropic exponential length space.  
If there is an
exponential curve with a conjugate point at $t_0$ 
then every
exponential curve from a point in $Y'$ has a conjugate point at $t_0$.
In fact so do all exponential curves in any covering space
of $Y$.
\end{lemma}

Note that this is not true for cut points as can be seen when
$Y$ is a cylinder.

\Pf
Fix $x,y \in Y$ and $v,\bar{v} \in S^{n-1}$.

If an exponential curve $exp_x(tv)$ has a conjugate 
point at $t_0$ then there exist distinct exponential curves
$exp_x(tv_i)$, $exp_x(tw_i)$ running from
$x$ to $exp_x(t_iv_i)=exp_x(s_iw_i)$ converging to it with $t_i, s_i \to t_0$
and $v_i, w_i \to v$.  
Eventually 
\be
\max \{d_S(v_i,v), d_S(w_i,v), d_S(v_i,w_i)\} < \min\{\theta_{(N+1)R/4,x,v, R/4}, \theta_{(N+1)R/4,y,w, R/4}\}
\ee
where we take $N > 8t_0/R$.  

Choose $\bar{v}_i$ and $\bar{w}_i$ in $S^{n-1}$ such that
\be
d_S(\bar{v},\bar{v}_i)=d_S({v},{v}_i)\qquad d_S(\bar{v},\bar{w}_i)=d_S({v},{w}_i)\qquad d_S(\bar{w}_i,\bar{v}_i)=d_S({w}_i,{v}_i)
\ee
Then by Lemma~\ref{ctogether} we have
\begin{eqnarray} \label{ctog2}
d_Y(exp_{x}(tv_i), exp_{x}(sv)) & = & d_Y(exp_{y}(t\bar{v}_i), exp_{y}(s\bar{v})) \\  %extra
d_Y(exp_{x}(tw_i), exp_{x}(sv)) & = & d_Y(exp_{y}(t\bar{w}_i), exp_{y}(s\bar{v})) \\  %extra
d_Y(exp_{x}(tv_i), exp_{x}(sw_i)) & = & d_Y(exp_{y}(t\bar{v}_i), exp_{y}(s\bar{w}_i)) 
\end{eqnarray}
whenever $t,s \in [0,NR/4]$ with $|t-s|<R/4$.

Since $t_i,s_i \to t_0$, eventually they are close and $\in [0, NR/4]$ thus
\be
0=d_Y(exp_{x}(t_iv_i), exp_{x}(s_iw_i)) =  d_Y(exp_{y}(t_i\bar{v}_i), exp_{y}(s_i\bar{w}_i)). 
\ee

\qed

%------------------------------------------------------------------------

\sect{Simply Connected Locally Isotropic Spaces} \label{simpconn}

In this section we study simply connected locally isotropic exponential length spaces
proving that they are homeomorphic to either $\RR^n$ or $\SSS^n$ 
[Theorem~\ref{homeomdone}]
and constructing global isometries [Lemmas~\ref{globisom} and \ref{isomext}].  We also prove a nice result
about triangles [Lemma~\ref{triangle}].

\begin{theorem} \label{homeomdone}
A simply connected locally isotropic exponential length space $Y$
is homeomorphic to $\RR^n$ or $\SSS^n$ where $n$ is the exponential
dimension of $Y$.
The homeomorphism is $exp_y:U \to Y$ where 
$U=\RR^n$ if $Y$ is unbounded and, if $Y$ is bounded,  
$U=B_0(D) \setminus \sim $ 
where $v \sim w$ if $|v|=|w|=D$  where 
$diam(Y)=D$.

\end{theorem}

\Pf
Clearly if there are no conjugate points then we are done
by Lemma~\ref{conjhomeom}.  Otherwise there is a first
conjugate point at some $t_0>R$ and by Lemma~\ref{allconj}
every exponential curve has a first conjugate point at the same $t_0>R$.
Choose one point $p\in Y$.  Then $\sup_{q\in Y} d_Y(p,q) \le t_0$
by Lemmas~\ref{minconj} and~\ref{mincut}. 

Since this is true for all $p\in Y$, we have
$t_0=diam(Y)$.  Furthermore $exp_p:B_0(t_0)\to B_p(t_0)$
is a homeomorphism by Lemma~\ref{conjhomeom}.  If $exp_p$ maps
$\partial B_0(t_0)$ to a single point then we are done.

Fix $v\in S^{n-1}$.  We know there exists $v_i \neq \bar{v}_i$
both converging to $v$ and $s_i, t_i \to t_0$ such that
$exp_p(t_iv_i)=exp_p(s_i\bar{v}_i)$.  Note that either
$t_i\ge t_0$ or $s_i\ge t_0$ since $exp_p$ is $1:1$ on
$B_0(t_0)$.  

Let $\theta_i=d_S(v_i,v)$
and $\bar{theta}_i=d_S(\bar{v}_i,v)$ so $\theta_i \to 0$ and 
$\bar{\theta_i} \to 0$.
Eventually 
\be
\max\{\theta_i, \bar{\theta}_i\} <\theta_{2t_0, p, R/4}/2.
\ee  

So applying Lemma~\ref{ctogether} taking $p$ to $p$, $v_i$ to itself and $\bar{v}_i$ to
some vector $w_i$,
we get $exp_p(t_iw_i)=exp_p(s_i\bar{v}_i)$.  So in fact
\be
exp_p(t_iw_i)=exp_p(t_iv_i).
\ee

Since $d_s(v_i,w_i)=\bar{\theta_i}\le \theta_i+\bar{\theta}_i \to 0$. 
we know eventually it is $< \theta_{2t_0, p,  R/4}/2$
Given any $\bar{v}, \bar{w} \in S^{n-1}$ such that $d_S(\bar{w},\bar{v})=\bar{\theta}_i$ we can
apply Lemma~\ref{ctogether} again mapping $p$ to $p$, $v_i$ to $\bar{v}$ and
$w_i$ to $\bar{w}$ to get 
\be \label{fixangle}
exp_p(t_i\bar{v})=exp_p(t_i\bar{w}) \textrm{ whenever } d_S(\bar{v}, \bar{w})<\bar{\theta}_i.
\ee

Now for any $w\in S^{n-1}$, $d_S(w,v_0)=k_i\bar{\theta}_i + \phi_i$ 
where $\phi_i<\theta_i$.  So there exists $\bar{w}_i \in S^{n-1}$
such that $d_S(w,w_i)<\theta_i$.  Applying (\ref{fixangle}) repeatedly along an exponential curve 
from $w_i$ to $v_0$ at intervals of length $\theta_i$, we get
$exp_p(t_iv)=exp_p(t_i\bar{w}_i)$.  Taking $i$ to infinity
and using the continuity of $exp_p$, $t_i\to t_0$, $w_i \to w$
we get $exp_y(t_0v_0)=exp_y(t_0w)$.
\qed

\begin{lemma} \label{globisom} 
Suppose $Y$ is a simply connected locally isotropic exponential length space.  
If $y_1,y_2\in Y$ and $g\in S0(n)$ then there is an isometry
$f_g:Y \to Y$ such that
$f_g(x)=exp_{y_2}(g(exp_{y_1}^{-1}(x)))$ for all $x\in B_y(D)$ where $D=diam(Y)\in (0,\infty]$.
\end{lemma}

\Pf 
By Lemma~\ref{homeomdone} we know that $exp_{y_1}^{-1}$ is well defined on
$B_{y_1}(t_0)$ where $t_0$  is both the first conjugate point and the diameter.
Here $t_0$ may be infinity.  We need only verify that
$f_g$ is an isometry from $B_{y_1}(D)\to B_{y_2}(D)$ since then it is forced
to be a global isometry by continuity since $\partial B_{y}(D)$ is a single point
for all $y\in Y$. 

Since $Y$ is simply connected we need only verify
that $f_g$ is a local isometry.  That is, for all $x\in Y$, there exists $r>0$
such that $f_g$ restricted to $B_x(r)$ maps isometrically onto $B_{f(x)}(r)$. 
We will choose $r< D-d_Y(x,y_1)$ to avoid trouble in the $\SSS^n$ case.

Now fix $x_1\in B_{y_1}(D)$ and let $s_1v_1=exp_{y_1}^{-1}(x_1)$.  Let 
\be
\theta_1 =\min\{\theta_{(D,y_1, R/4}, \theta_{D,y_2, R/4}\}, 
\ee
Since $exp_{y_1}$ is a homeomorphism, we can take 
$r$ to be sufficiently small that for all $x\in B_{x_1}(r)$, 
we have $d_S(exp_{y_1}^{-1}(x)/|exp_{y_1}^{-1}(x)|, v_1)<\theta_1$.   

For any $z_1, z_2\in B_{x_1}(r)$, let $s_iw_i=exp_{y_1}^{-1}(z_i)$.
Since $d_S(w_i,v_1)<\theta$, we can apply Lemma~\ref{ctogether}
to see that 
\begin{eqnarray}
d_Y(z_1,z_2)&=&d_Y(exp_{y_1}(s_1w_1),exp_{y_1}(s_2w_2))\\
&=&d_Y(exp_{y_2}(s_1gw_1), exp_y(s_2gw_2))=d_Y(f_g(z_1), f_g(z_2)).
\end{eqnarray}
\qed

Note that the above Lemma implies that if we have a triangle formed by two length
minimizing curves of lengths $a<D$ and $b<D$ and any angle $\theta$ between them then the
length of the third side is determined depending only on $\theta$, $a$ and $b$.

\begin{coro} \label{Disot}
Suppose $Y$ is a simply connected locally isotropic exponential length space
then its isotropy radius is the diameter $D$ or infinity in the unbounded case.
\end{coro}

We now prove that given a triangle
with sides of length $a$, $b$ and $c$ we can determine the angle opposite $c$,
however we must restrict our lengths to avoid the pole which causes 
indeterminacy
even in $\SSS^n$.   
 
\begin{lemma} \label{triangle}
If $Y$ is a simply connected
locally isotropic exponential length space with diameter $D\in [R, \infty]$ then for all 
$y\in Y$, $a,b,c\in (0,D)$, if
$d_Y(exp_y(bv), exp_y(aw)) = c$ then $d_S(v,w)=\theta(a,b,c)$.
In particular, if $a+b+c=2D$ then $\theta(a,b,c)=\pi$.
\end{lemma}

\Pf
Corollary~\ref{Disot} and Lemma~\ref{solveF} imply the existence of $\theta(a,b,c)$.
%before december I redid the whole proof

Look at the triangle between the points $x_2$, $exp_{x_2}(a{w}_2)$
and $exp_{x_2}(bv_2)$.  Join the later two points by a length minimizing
curve $C(t)$ parametrized by arclength
such that $C(0)=exp_{x_2}(bv_2)$.  Note that since $d_Y(x_2, C(t))\le D$ then
$a+t\le D$ and $b+(c-t)\le D$ so $a+b+c \le 2D$.

If $a+b+c=2D$ then  $C(t)$ hits the point $\bar{x}_2=\partial B_{x_2}(D)$.
Since $C(t)$ runs minimally to this point and $d_(exp_{x_2}(tw_2), \bar{x}_2)=D-t$ including $t=a$,
we know $C(t)=exp_{x_2}((t+a)w_2)$ for $t\in [0,D-a]$.  Similarly
$C(c-t)=exp_{x_2}((t+b)v_2)$.  So we have an exponential curve running from $x_2$ through $\bar{x}_2$ and
back to $x_2$ and $a+b+c=2D$.  Since this curve must be minimizing on segments of length $D=t_0$ by the
by Corollary~\ref{conjmin} we have $d_Y(exp_{x_2}((D/2)v_2), exp_{x_2}((D/2)w_2))=D$. However, we can 
join these points by a curve through $x_2$ of length $D$ so that curve must be an exponential curve
and so $d_S(v_2, w_2)=\pi$ by Lemma~\ref{pi}.

\qed

We can now use Lemma~\ref{ballext} to extend isometries between subdomains of $Y$
to all of $Y$.

\begin{theorem}  \label{isomext}
If $Y$ is a locally isotropic exponential length space
which is simply connected.  If $f:A\to B$ is an isometry
between subsets of $Y$, then there is an extension of $f$
to an isometry from $Y$ to $Y$.
\end{theorem}

\Pf
If $Y$ is unbounded then by Corollary~\ref{Disot}, $R=\infty$ and
this is a consequence of Lemma~\ref{ballext}.  Suppose $Y$ is bounded.

If $A=Y$ then there is nothing to do.  If there exists $x\in Y\setminus A$
then $A \subset B_y(D)$ where $d_Y(y,x)=D$.  So we can apply Corollary~\ref{Disot}
and Lemma~\ref{ballext}, to extend $f$ to $B_y(D)$.  Let $f(x)$ be the
only point in $Y\setminus B_x(D)$.  Then $f$ is continuous so it must be an isometry.
\qed

As a consequence of this theorem we know that there are global isometries
mapping any point to any other point, that there are global isometries
mapping any triangle to any other congruent triangle.  In the unbounded case, combined
with the fact that there the space is globally minimizing, we can use Busemann's Theorem
to state that $Y$ is either Euclidean or Hyperbolic space \cite{Bu}.

%-----------------------------------------------------------------------------

%could circumvent Birkhoff perhaps
\sect{Birkhoff's Theorem and Local Isotropy} \label{sectbirk}

In this section we apply Birkhoff's Theorem to complete proofs
of Theorems~\ref{isoexplthm} and~{\ref{MainThm2} \cite{Birk}.
This theorem dates to 1941 and there are similar earlier theorems
by Busemann which characterize $\HH^n$ and $\EE^n$ \cite{Bu}.  The difficulty
with Busemann's theorems is that they assume extendibility of the geodesics
as minimizing curves and this is not true in the sphere.  
Birkhoff's Theorem is stated in the following proof.

\begin{theorem}\label{birkhoff}
If $Y$ is a locally isotropic exponential length space
which is simply connected then $Y$ is isometric to $\SSS^n$,
$\HH^n$ or $\RR^n$ where $n$ is the exponential dimension of $Y$.
\end{theorem}

\Pf
By Birkhoff's Theorem \cite{Birk} any length space such that
the following hold is such a simply connected space form:

{\em
a) $\bar{Y}$ has locally unique minimal geodesics. 

b) any isometry on subsets of $\bar{Y}$ extends to an isometry of the whole
space. }

\noindent
Now (a) follows from Lemma~\ref{localmin}  %MORE DETAILS?
and (b) follows from Theorem~\ref{isomext}.  Then we apply Theorem~\ref{homeomdone}
to show that $n$ matches the dimension of the space form.
\qed

\begin{remark}
We could also prove this theorem without quoting Birkhoff, but rather using a more
recent theorem which states that all locally compact two point homogeneous manifolds 
are Riemannian manifolds.  
Our space is locally compact by the definition of exponential length space and
it can be seen to be two point homogeneous by taking $A$ to contain exactly
two points in Theorem~\ref{isomext}.  Once $Y$ is isometric to a Riemannian
manifold, we know the exponential curves are geodesics because they are locally
minimizing.  We also know it is an $n$ dimensional manifold by Theorem~\ref{homeomdone}
where $n$ is the exponential dimension of $Y$.
Then we can apply Theorem~\ref{isomext} again with $A$ and $B$ sharing
a geodesic and forcing the existence of isometries rotating around that geodesic
to get constant sectional curvature.  All simply connected space forms must
be either $\SSS^n$, $\HH^n$ or $\EE^n$ (c.f. \cite{doC}).
\end{remark}

This would appear to be sufficient to complete the paper but we must relate the
exponential structure of the space form to that of the exponential length
space.  

\begin{lemma} \label{match}
If $X$ and $Y$ are locally isotropic exponential length spaces
off $W_X$ and $W_Y$ that are  isometric to each other
then they have the same exponential length structure.
That is $f:Y \to X$ an isometry implies that for all $y \in Y\setminus W_Y$ 
mapped to $f(y)\in X \setminus W_X$ there exists 
$g_y\in S0(n)$ such that 
\be
f(exp_y(v))=exp_{f(x)}(g_y v).
\ee
and if $R$ is the minimum of the isotropy radii of $X$ and $Y$, we have
\be 
F_{y}(\theta,s,t)=F_{f(y)}(\theta,s,t) \qquad \forall s,t<R.
\ee
\end{lemma}

Note that without assuming local isotropy on $X$ this is false since
$\RR^n$ can be given two distinct exponential length structures
that are both isometric to Euclidean space.  See Example~\ref{nomatch}.
Note also that $g_x$ need not be continuous in $x$.

\Pf  %old proof using M before december
First since we are only making a statement about $y\in Y\setminus W_Y$ and
$f(y)\in X \setminus W_X$ we can first choose $Y'$ and $X'$ to be their
respective connected components and then extend the exponential
structures to $Cl(Y')$ and $Cl(X')$ respectively using Lemma~\ref{badsexp}.

Thus, without loss of generality, we may assume $W_Y$ and $W_X$ are
empty sets and can then apply all our lemmas concerning such spaces,
using the fact that they are extended exponential length spaces.

Since $\gamma(t)=exp_y(tv)$ is a length minimizing curve for $t\in (0,R)$, 
$f(\gamma(t))$ must be as well.  So $f(\gamma(t))$ is an exponential curve starting at
$f(p)$.  Thus there is a map $g_p:\RR^n \to \RR^n$
such that $f(exp_p(v))=exp_{f(p)}(g_p(v))$.  
By Lemma~\ref{pi} we have
$g_p(tv)=tg_p(v)$ even for negative $t$.  
Using the fact that exponential curves are parametrized proportional
to arclength 
we see that $|g_p(v)|=L(f(\gamma([0,L]))=L(\gamma([0,L]))=|v|$.
So we need only verify that $g_p$ is an isometry from $S^{n-1}$ to $S^{n-1}$.

By Lemma~\ref{pi} we know that for $\epsilon>0$ sufficiently small
\be
F_y(\pi, \epsilon, \epsilon)=d_Y(exp_y(\epsilon v), exp_y(-\epsilon v))=2 \epsilon
\ee
and
\be
F_{f(y}(\pi, \epsilon, \epsilon)=d_Y(exp_{f(y)}(\epsilon v), exp_{f(y)}(-\epsilon v))=2 \epsilon.
\ee

We will prove that for any isometry $f$ and $g_f$ as above,
\be \label{induct1}
F_y(\pi/2^k, \epsilon, \epsilon)=F_{f(y}(\pi/2^k, \epsilon, \epsilon)
\ee
and
\be \label{induct2}
d_S(v_1, v_2)=\pi/2^k \textrm{ implies that } d_S(g_f(v_1), g_f(v_2))=\pi/2^k
\ee 
by induction on $k$.  

When $k=0$, (\ref{induct1}) holds by Lemma~\ref{pi} as described above.  
On the other hand $d_S(v_1,v_2)=\pi$ implies $v_1=-v_2$ so $g_f(v_1)=-g_f(v_2)$
and we are done.

Assuming it is true for $k=j$, let $v_1, v_2 \in S^{n-1}$
be chosen such that $d_S(v_1,v_2)=\pi/2^j$.  Let $w=(v_1+v_2)/|v_1+v_2|$.
So
\be
d_Y(exp_y(\epsilon v_1), exp_y(\epsilon w))=F_y(\pi/2^{j+1})
d_Y(exp_y(\epsilon v_2), exp_y(\epsilon w))
\ee
Thus by the isometry,
\be
d_X(exp_{f(y)}(\epsilon g_f(v_1)), exp_{f(y)}(\epsilon g_f(w)))=F_y(\pi/2^{j+1})
=d_X(exp_{f(y)}(\epsilon g_f(v_2)), exp_{f(y)}(\epsilon g_f(w))).
\ee
Thus by Lemma~\ref{solveF}, $d_S(g_f(v_1),g_f(w))=d_S(g_f(v_2),g_f(w))$.
Since the triangle inequality gives, $d_S(g_fv_1,g_fw)+d_S(g_fv_2,g_fw) \ge d_S(g_fv_1, g_fv_2)$,
we know $d_S(g_fv_1,g_fw)=d_S(g_fv_2,g_fw) \ge \pi/2^{j+1}$.  By
the properties of $F_{f(y)}$ this implies that
\begin{eqnarray}
F_{f(y)}(\pi/2^{j+1}, \epsilon, \epsilon)
&\le& F_{f(y)}(d_S(g_fv_1, g_fw), \epsilon, \epsilon) \\
&=& d_X(exp_{f(y)}(\epsilon g_f(v_1)), exp_{f(y)}(\epsilon g_f(w)))
=F_y(\pi/2^{j+1}).
\end{eqnarray}
However, the same holds for the isometry $f^{-1}$ so we get the opposite inequality.
Thus we get (\ref{induct1}) for $k=j+1$ and then Lemma~\ref{solveF} gives us (\ref{induct2}).

Thus applying (\ref{induct2}) and the properties of $g_f$ shown at the top, we get,
\begin{eqnarray} \label{induct2b}
F_y(\pi/2^k, s,t)& =& d_Y(exp_y(tv_1), exp_y(sv_2))=d_X(exp_{f(y)}(tg_fv_1), exp_{f(y)}(s g_fv_2)) \\
&= & F_{f(y)}(\pi/2^k, s,t) \qquad \forall s,t < \min\{R_y, R_{f(y)}\}.
\end{eqnarray}

Now choose any $j\in \NN$ such that $j\pi/2^k<\pi$
and any $v_1, v_2$ with $d_S(v_1, v_2)=j\pi/2^k <\pi$.  Let $c$ be a length minimizing
curve from $c(0)=exp_y(tv_1)$ to $c(L)=exp_y(sv_2)$ where $L=F_y(j\pi/2^k, s,t)$.
Note that $L<R$, so $c$ cannot leave $B_y(R)$. 
Then $f(c(t))$ is length minimizing from $exp_{f(y)}(tg_fv_1)$ to $exp_{f(y)}(sg_fv_2)$ 
so $L=F_{f(y)}(d_S(g_fv_1, g_fv_2), s,t)$.  

We claim that $v_t=exp_y^{-1}(c(t))/|exp_y^{-1}(c(t))|$  is in the minimizing geodesic segment between
$v_1$ and $v_2$ for all $t\in [0,L]$.  If not at some point $t$,  there is an element $g\in S0(n)$
which maps $v_i$ to $v_i$ but moves $v_t$.  Then 
the isometry $f_g$
of Lemma~\ref{globisom} maps $c(t)$ to another length minimizing curve with the same end points.
But our space is locally minimizing, so there cannot be a second such curve and we have a contradiction.

Thus we can choose $t_0=0<t_1< t_2<...< t_j=L$
so that $c(t_h)=exp_y(s_h w_h)$ where  $d_S(w_h, w_{h+1})=\pi/2^k$.
This can be done using the
intermediate value theorem and continuity of $exp_y^{-1}$ on $B_y(R)$.  
By (\ref{induct2})
$d_S(g_fw_h, g_f w_{h+1})=\pi/2^k$, so in fact
\be \to 
d_S(g_f v_1, g_f v_2) \le \sum_{i=0}^{j-1} d_S(g_fw_h, g_f w_{h+1})=j\pi/2^k.
\ee
Since this whole argument works for $f^{-1}$ as well we get
\be \label{induct4}
d_S(v_1, v_2)=j\pi/2^k \textrm{ implies that } d_S(g_rf(v_1), g_f(v_2))=j\pi/2^k.
\ee 
Using this fact and the fact that $L$ was preserved under the isometry $f$, we get
\be \label{induct3}
F_y(j\pi/2^k, s,t)=L=F_{f(y}(j\pi/2^k, s,t) \qquad \forall s,t < R=\min\{R_y, R_{f(y)}\}.
\ee

Now given any $\theta \in [0,\pi]$ and any $v, w$ such that $d_S(v,w)=\theta$ there
exists $v_i, w_i$ such that $d_S(v_i, w_i)= j_i\pi/2^{k_i}$ so that substituting these
$v_i$ and $w_i$ in (\ref{induct3}) and (\ref{induct4}) and taking $i\to \infty$ we
get 
\be \label{induct5}
d_S(v, w)= d_S(g_rf(v), g_f(w))
\ee 
and
\be \label{induct6}
F_y(\theta, s,t)=F_{f(y}(\theta, s,t) \qquad \forall s,t < R=\min\{R_y, R_{f(y)}\}.
\ee
\qed

We can now prove that a locally isotropic exponential length space is a collection of manifolds with constant
sectional curvature joined at discrete points.

\noindent{\bf Proof of Theorem~\ref{isoexplthm}:}

For the first part we
need only show that if $\bar{Y}$ is the universal cover of $(Cl(Y'))$ 
then it is either $S^n$, $H^n$ or $E^n$ with
their standard metrics.  By Lemma~\ref{badsexp}, we 
know the $Cl(Y')$ is isometric to a locally isometric exponential length
space and so by Lemma~\ref{univcov} its universal cover $\bar{Y}$ exists and
is a simply connected locally
isometric exponential length space.  So by Theorem~\ref{birkhoff}, $\bar{Y}$ is 
isometric to $\SSS^n, \HH^n$ or $\EE^n$.  Thus $Cl(Y')$ is isometric to
a space form of some constant curvature $K$.  Lemma~\ref{match} then says that  
the exponential structures match and
$F_y(\theta,s,t)=F_K(\theta,s,t)$ for all $y \in Y'$.

Once this is known we use the fact that $W_Y$ is discrete, to
piece together the various connected components of $Y'\cap W_Y$
in a countable way.
\qed

Finally we can prove Theorem~\ref{MainThm2} which implies Theorem~\ref{MainThm}.

\noindent{\bf Proof of Theorem~\ref{MainThm2}:}  
By Gromov's Compactness Theorem and the $f$ ball packing property,
we know that a subsequence of these $M_i$ (also called $M_i$)
converges to some complete length space $Y$.  

By Theorem~\ref{allbutonto}, Lemma~\ref{discrete}
and Lemma~\ref{ontoconn}, we know
that $Y$ is a locally isotropic exponential length space off a discrete
set $W_Y$.    We can thus apply Theorem~\ref{isoexplthm} to 
obtain the required properties of $Y$.

Since the exponential length structure originally given to $Y$ as a limit
of the $M_i$ was defined to satisfy
\be
\lim_{i\to\infty}F_{q_i}(\theta, s,t)=F_y(\theta,s,t)
\ee
and since Theorem~\ref{isoexplthm} says that the exponential length structure
on any $Cl(Y')$ must match that of a space form with constant sectional curvature
$K$, we see that $F_y(\theta,s,t)=F_{K}(\theta,s,t)$.  This implies (\ref{MainThm2a}). 
\qed

%-------------

\sect{Ricci Curvature} \label{sectricci}
In this section we will apply Theorem~\ref{MainThm2} combined with the lower Ricci curvature 
bound to prove Theorem~\ref{RicciThm}.

\noindent{\bf Proof of~\ref{RicciThm}:}
First note that if $M$ has a fixed lower bound on Ricci curvature
then by Bishop-Gromov, it satisfies the conditions of Theorem~\ref{MainThm}.
Furthermore $Y$ can be described as a limit of $M_i$ satisfying this
uniform Ricci curvature bound just as in Theorem~\ref{MainThm2}
except that now we have additional measure properties on $Y$
proven by Colding.  We will show $Y$ is a space form by showing $Y=Cl(Y')$.

If $Y\neq Cl(Y')$ then by Theorem~\ref{MainThm} we know that 
$Y$ contains at least two
space forms $Cl(Y')$ and $Cl(Y'')$ joined at a common point $y_0$.
Let $p_i \to y$.  We know $vol(B_{p_i}(r))$ converges to $vol(B_y(r))$
by Colding's Volume Convergence Theorem \cite{Co}.  
By the properties of $Cl(Y')$ and $Cl(Y'')$, we know 
\be
vol(B_y(r))\ge V(n,K',r)+V(n,K'',r) \forall r>0.
\ee
On the other hand, since we don't have dense bad points then there exists
$y_1$ near $y$ in $Y'$, such that $vol(B_{y_1}(r))=V(n,K',r)$ for all $r$
sufficiently small.  
Let $d(y_1,y_0)=r_1$.
Take an annulus about $y_1$ which includes $y_0$ but no other bad points.
By volume comparison with $H=\min\{K', K''\}$:
\be
Vol(Ann_{y_1}(r_1-r, r_1+r))\le 
\frac{V(n,H,r_1+r)-V(n,H,r_1-r)}{V(n,H,r_1-r)} Vol(B_x(r_1-r)).
\ee
By the space forms:
\begin{eqnarray}
Vol(Ann_{y_1}(r_1-r, r_1+r))&\ge& Vol(Ann\cap Y')+Vol(Ann\cap Y'')\\
&\ge& (V(n,K',r_1+r)-V(n,K',r_1-r)) + V(n,K'',r)
\end{eqnarray}
Putting this together and using $r_1-r<r_x$, 
$
Vol(B_x(r_1-r))=V(n,K',r_1-r)
$
we have
\be
\frac{V(n,H,r_1+r)-V(n,H,r_1-r)}{V(n,H,r_1-r)} 
\ge \frac{V(n,K',r_1+r)-V(n,K',r_1-r)) + V(n,K'',r)}{V(n,K',r_1-r) }
\ee
Now we can take any $y_1$ close to $y_0$ and set $r_1=2r$.
\be
\frac{V(n,H, 3r)-V(n,H,r)}{V(n,H,r)} 
\ge \frac{V(n,K',3r)-V(n,K',r)) + V(n,K'',r)}{V(n,K',r) }
\ee
Taking $r\to 0$ and using the fact $\lim_{r\to 0}V(n, H, r)/r^n=\omega_n$,
we get the impossible limit:
\be
\frac{3^n-1^n}{1^n} 
\ge \frac{3^n-1^n + 1^n}{1^n}.
\ee
\qed

%-----------------------------

\end{document}